\batchmode
\documentclass[12pt,titlepage]{article}
\usepackage{amsmath, amssymb, amsthm, amsfonts}
\numberwithin{equation}{section}
\newtheorem{thm}{Theorem}[section] 
\newtheorem{theorem}[thm]{Theorem} 
\newtheorem{definition}[thm]{Definition}
\newtheorem{df}[thm]{Definition}
\newtheorem{dfn}[thm]{Definition}
\newtheorem{defn}[thm]{Definition}
\newtheorem{lemma}[thm]{Lemma}
\newtheorem{lem}[thm]{Lemma}
\newtheorem{lm}[thm]{Lemma}
{\newtheorem{remark}[thm]{Remark}}
{\newtheorem{rem}[thm]{Remark}} 
{\newtheorem{example}[thm]{Example}}
{\newtheorem{ex}[thm]{Example}}
\newtheorem{cor}[thm]{Corollary}
\newtheorem{prop}[thm]{Proposition}
\def\pf{{\em Proof}.\, }
\def\qed{\hfill $\square$\par\vspace{5pt}}
\def\bysame{\leavevmode\hbox to3em{\hrulefill}\,}
\def\a{\mathfrak{a}}

\def\g{\mathfrak{g}} 
\def\h{\mathfrak{h}}

\def\t{\mathfrak{t}}

\def\N{{\bf N}}
\def\R{{\bf R}}
\def\C{{\bf C}}

\def\H{{\bf H}}
\def\Z{{\bf Z}}
\def\F{\mathcal{F}}
\def\J{\mathcal{J}}
\def\L{\mathcal{L}}
\def\i{\sqrt{-1}}
\renewcommand{\Re}{\mbox{\rm Re}\,}

\title{Lectures on Dunkl operators}  
\author{Eric M. Opdam 
\thanks{These lecture notes are based on a series of lectures given by 
E.~Opdam in  the project research ``Harmonic analysis on homogeneous 
spaces and representation of Lie groups'' 
at RIMS, Kyoto University (Japan) in 1997.  
He gave five lectures
from 27 to 30 October and 25 November 1997 on trigonometric 
Dunkl operators, degenerate affine Hecke algebra, and harmonic
analysis for the hypergeometric function for root systems.  
These notes are prepared  by T.~Honda (section 3, 4, 5), 
H.~Ochiai (section 2, 6), N.~Shimeno (section 8, 9), and K.~Taniguchi 
(section 7) after Opdam's lectures. }
}
\date{}
\makeindex
\begin{document}
\maketitle
\tableofcontents
\newpage
\section{Prologue}
These are the lecture notes of a series of 5 lectures 
held at RIMS in October/November 1997, 
in one of the workshops of the research project 
``Harmonic analysis on homogeneous spaces and representation of 
Lie groups''. In these lectures I have discussed Dunkl operators 
in the trigonometric, differential setting. This 
subject has been very dear to me for many years, and  
it was a great pleasure to have the opportunity    
to lecture on this subject in a stimulating environment.
My warm thanks go out to those who made this possible: to prof. T. Oshima  
for inviting me to participate in the research project 
``Harmonic analysis on homogeneous spaces and representation of 
Lie groups'' at RIMS; to prof. M. Kashiwara 
for being my host at the RIMS institute; and to the note takers  
T.~Honda, H.~Ochiai, N.~Shimeno, and K.~Taniguchi for their 
kindness to prepare these notes.
 
The choice of the subject is based on my personal 
experience and taste. In view of the recent developments 
concerning Dunkl operators,   
one may object that my choice represents  
a rather limited point of view. 
Indeed, in view of Cherednik's work, the trigonometric 
differential Dunkl operators seem  
to be only a degenerate limit of a theory of commuting 
difference operators. These difference operators arise from 
commutation formulae inside Cherednik's double affine Hecke 
algebra. This magnificent insight has changed the way in which 
we ought to think about Dunkl operators and their applications.

Nonetheless, I have restricted myself to discuss the 
differential case. There are various reasons for doing so. 
First of all, there are a number of recent expositions 
(\cite{Kir}, \cite{Mac}, \cite{C2}) 
of the new algebraic 
theory of Dunkl operators and the double affine Hecke algebra.
Second, 
the trigonometric differential limit that we consider, is very rich  
and it has served as a guideline 
for developments in the general theory. Third, there are 
aspects in the differential theory that have resisted generalization 
to the general theory so far. 
Especially with respect to harmonic analysis, the 
differential theory has currently reached a higher level of maturity 
(although an exciting start of the harmonic analysis for the 
difference equations can be found in \cite{C3}).
It is this analytic aspect of the theory of Dunkl operators I shall 
concentrate on. Finally, although we will only deal with the differential 
theory, on our way we shall meet with the (degenerated) double 
affine Hecke algebra several times. As has been mentioned before, Cherednik's 
approach has profoundly changed our perception of Dunkl operators, and 
of course this also manifests itself in the differential theory. In 
fact, I hope and even expect that for some readers, this modern treatment of 
Dunkl operators 
will be a motivation to look more closely at the double affine Hecke 
algebra.

Let me make some personal historical comments on the development 
of the theory we will be studying in these notes.  
Dunkl operators were conceived by Charles
Dunkl in 1989 (see \cite{Duone}).  He found  these operators in the so-called
rational differential situation,  which is the  basic example. He proved the
two fundamental properties,  the $W$-equivariance (which is in fact immediate
here) and the marvellous  commutativity, and he used this to set up a theory
analogous to the theory  of spherical harmonics. 

Almost at the same time, but unaware of Dunkl's 
fundamental results, 
Gerrit Heckman and I were seeking to generalize the theory of 
the spherical function of Harish-Chandra. Our goal was a theory of 
multivariable hypergeometric functions associated with a root system. 
Inspired 
by Tom Koornwinder's work \cite{K} in this direction (already in the 
early seventies) we set up such 
a theory in a series of papers \cite{HO1}, \cite{Hec1}, \cite{Op0},
\cite{Op1}. 

Soon afterwards I noticed (\cite{Op02}) that this theory provided
natural tools (shift operators) that could be succesfully applied to a number of
combinatorial and analytic problems that were related to root systems 
(most notably, Macdonald's constant term conjectures for root systems
\cite{Matwo}).  In spite of these applications,  the hypergeometric theory 
itself was not in a very satisfactory state
at the time.  The main arguments were indirect and complicated, avoiding
at all times to use explicit knowledge of the defining differential
equations of our  hypergeometric function.
The obstacle, psychologically,  was that  it {\it seemed} hopeless to write
down these defining differential equations explicitly, 
since this was already  impossible (in general) for
Harish-Chandra's spherical function itself. 

These difficulties were resolved
in a rather drastic way when Gerrit Heckman noticed (\cite{Hec2.5})  the
connection with Dunkl's work. Dunkl's operators  provided a very simple method
for constructing the differential equations  we needed, in the rational
version of our theory. Heckman defined a trigonometric  version of these
operators as well (\cite{Hec3}). There was
however a remarkable difference with the  rational case: the trigonometric 
operators that
Heckman found were $W$-equivariant,  but they did not commute. Nonetheless
these ``Dunkl-Heckman'' operators were  important and useful, because they were
the building blocks for the desired  commuting (higher order) differential
operators (and shift operators)  in the trigonometric case. 

The next development was Ivan Cherednik's discovery of the connection 
between (degenerated) affine Hecke algebras on the one hand,   
and Dunkl and Dunkl-Heckman operators on the other hand (\cite{Chone}, 
\cite{Chtwo}). 
This discovery  had some important consequences. From the structure 
theory of 
Hecke algebras it was now obvious that there also existed 
{\it commuting} 
Dunkl-type operators in the trigonometric case. 
It is an interesting fact that these commuting operators are not  
W-equivariant in the trigonometric case. 
The joint spectral theory of these commuting 
``Dunkl-Cherednik'' operators will be the main subject of study 
in these notes. Noncompact spectral theory started with  
De Jeu's important paper \cite{J} (the rational case), and was then 
further explored in the 
trigonometric case in \cite{Op}, \cite{Op2} and in
Cherednik's paper \cite{C1}.

Cherednik's discovery also created a natural way to discretize the 
theory (creating the difference operators alluded to in the 
second paragraph of this prologue),  by using the affine Hecke algebra 
instead of the degenerated version. This led to the complete 
solution of the Macdonald conjecures (including the ``q-version''), and 
many new results (see \cite{C2} for a very good account of these
developments). 

\section{Dunkl operators in the trigonometric setting}
The basic reference for this section is \cite{Op}.

\subsection{Notation}
We assume that the reader is familiar with root 
systems and their basic properties. However, in order to fix notations and
conventions we will review the definitions of these and related 
fundamental structures in this subsection.

Let \index{a@$\mathfrak a$} $\mathfrak a$ be a  Euclidean vector space of
dimension $n$. For $\alpha\in \mathfrak a^*$ we denote by
$X_\alpha\in \mathfrak{a}$ the element corresponding to $\alpha$.
When $\alpha$ is nonzero we introduce the covector $\alpha^\vee\in\mathfrak a$ 
of $\alpha$ by the formula
\[
\alpha^\vee =\frac{2 X_\alpha}{ ( X_\alpha, X_\alpha ) }.
\]

A nonzero $\alpha$ in $\mathfrak a^*$ determines the 
orthogonal reflection \index{ra@$r_\alpha$} $r_\alpha \in O(\mathfrak{a})$
in the hyperplane $\mathrm{ker}(\alpha)$ of $\mathfrak a$. This reflection 
is given by the formula
\[
r_\alpha(\xi)=\xi - \alpha(\xi)\alpha^\vee.
\]
In many instances the orthogonal transformation $r_\alpha$ will act on 
spaces derived from $\mathfrak a$, such as the complexification
of $\mathfrak a$, certain stable lattices in $\mathfrak a$, 
tori that are a quotient of $\mathfrak a$ by such a stable lattice, and 
also on the dual $\mathfrak a^*$. In all these situations we will simply use
the  same notation $r_\alpha$, since there is no danger of confusion (in the
last case, notice that $r_\alpha=r_{\alpha^\vee}$ when we identify
$\mathfrak a$ and $\mathfrak a^*$).

A finite subset \index{r@$R$}$R\subset {\mathfrak a}^*\backslash\{0\}$ is
called a root system when it satisfies the following properties:
\begin{itemize}
\item[(R1)] $R$ spans ${\mathfrak a}^*$.
\item[(R2)] $\forall\alpha\in R$,\  $r_\alpha(R)=R$.
\item[(R3)] $\forall\alpha, \beta\in R$,\  $\alpha(\beta^\vee)\in\Z$.
\end{itemize}

The elements of $R$ are called roots.
We shall always assume that $R$ is reduced%
\footnote{
This assumption is not necessary.
Actually, an important class of orthogonal polynomials (Koornwinder-
polynomials)
arises from the non-reduced root system of type BC$_n$.
However, we employ this assumption for simplicity.}, which means 
that $\R\alpha\cap R=\pm\alpha$ for every $\alpha\in R$.

Clearly the set \index{rv@$R^\vee$} $R^\vee=\{\alpha^\vee\mid\alpha\in
R\}\subset \mathfrak{a}$ is also a root system, called dual or coroot system.

The group generated by the reflections $r_\alpha$ is a finite reflection 
group, called the Weyl group and denoted by \index{w@$W$} \index{wr@$W(R)$}
$W=W(R)$. Because of (R3), 
\index{q@$Q$} \index{qa@$Q(R)$} $Q=Q(R)=\Z R$ and 
\index{q@$Q^\vee$} $Q^\vee =Q(R^\vee)$  \index{qrv@$Q(R^\vee)$} are stable lattices for
the action of $W$. 
These lattices are called the {\it root lattice} and the {\it coroot lattice} respectively. 
The dual lattice  \index{p@$P$} $P = \mbox{Hom}_{\bf
Z}(Q^\vee, {\bf Z}) \subset \mathfrak{a}^*$ is called the weight lattice of
$R$, and is of course also $W$ stable. 

We put \index{h@$\mathfrak{h}$}$
\mathfrak{h}= \mathfrak{a}_{{\bf C}}$ and   $\t=\sqrt{-1} \mathfrak{a}$, hence
we have $\h=\a+\t$.  \index{t@$ \mathfrak{t}$}
Let \index{h@$H$} $H$ be the complex torus $H = \mbox{Hom}_\Z(P,{\bf
C}^\times)  = Q^\vee \otimes_{{\bf Z}} {\bf C}^\times$. 
The Weyl group $W$ stabilizes $P$ and $Q^\vee$,
hence $W$ also acts on $H$. We have $H=TA$, 
where $T$ is a compact torus and $A$ is the real split 
torus, corresponding to $\t$ and $\a$ in $\h$ respectively. 
\index{a@$A$} \index{t@$T$} 

Choose and fix a halfspace in $\mathfrak a^*$ such that none of the roots  
of $R$ are in the boundary of this halfspace. The roots in this
halfspace are said to be positive, and the set of positive roots is
called a positive subsystem \index{rp@$R_+$} $R_+\subset R$. Let
\index{qp@$Q_+$} $Q_+$ be the ${\bf Z}_+$-span of $R_+$. It is well known that 
$Q_+$ is a simplicial cone over ${\bf Z}_+$, and is generated over ${\bf Z}_+$
by a basis of roots $\{ \alpha_1,\dots, \alpha_n \}$. Put $r_i = r_{\alpha_i}$,
then \index{s@$S$} $S = \{ r_1, \dots, r_n \}$ is 
a set of generators of $W$. In fact these generaters give a presentation of
$W$ as a Coxeter group, with relations $r_i^2=1$ and $(r_ir_j)^{m_{ij}}=1$.

The set $Q_+$ defines an important partial ordering $<$ in $\mathfrak a^*$ by 
$\lambda<\mu$ iff $\mu-\lambda\in Q_+$. This ordering is called the dominance
ordering. 

When $\lambda(\alpha_i^\vee)\geq 0 \forall i\in \{1,\dots,n\}$ we call
$\lambda$ dominant (and we call $\lambda$ strongly dominant when all the
inequalities are strict). The set $\mathfrak a^*_+$ of all strongly dominant
elements is  called the Weyl chamber. It is well known that the closure of the
Weyl chamber $\overline{\mathfrak a^*_+}$ is a fundamental domain for the
action of $W$. Let $P_+\subset P$ denote the set of  dominant
weights\index{pp@$P_+$}.  It is generated over $\Z_+$ by the basis
$\{\lambda_i\}$ dual to $\{\alpha_i^\vee\}$. The weights $\lambda_i$ are
called fundamental weights. 

Let \index{ch@${\bf C}[H]$} ${\bf C}[H]$ be the
space of Laurent polynomials (finite linear combinations of
algebraic characters $e^\lambda$ with $\lambda \in P$). By restiction to $T$
one may identify this space of functions with  the space of Fourier polynomial
on $T$. 

\subsection{Dunkl-Cherednik operator} 
\begin{prop}\label{prop-1.1}
The divided difference operator 
$\displaystyle\frac{1}{1-e^{-\alpha}}(1-r_\alpha)$ maps 
$\C[H]$ into itself. 
\end{prop}
\pf 
This easily follows from the summation over geometric series.
This operator sends 
\[
e^{\lambda} \mapsto
\left\{
\begin{array}{lll}
e^\lambda(1+e^{-\alpha}+\cdots+e^{(1-\lambda(\alpha^\vee))\alpha})
& \mbox{if} & \lambda(\alpha^\vee)>0 \\
0 & \mbox{if} & \lambda(\alpha^\vee)=0 \\
- e^{r_\alpha \lambda}
(1+e^{-\alpha}+\cdots+e^{(1+\lambda(\alpha^\vee))\alpha})
& \mbox{if} & \lambda(\alpha^\vee)<0 
\end{array}.
\right.
\]
This proves the required property.
\qed 

Notice the asymmetry, the difference between the formulae 
for positive exponents and for negative exponents. 
Only the largest element of $\lambda$ and $r_\alpha \lambda$,
(in the dominance order) shows up in the support of the image of 
$e^\lambda$. 
This property plays an important role in the sequel.

Let us introduce the Weyl denominator 
\index{delta@$\Delta$}
\[
\Delta=
{\displaystyle\prod_{\alpha\in R_+}}
\left( e^{\alpha/2} - e^{-\alpha/2} \right)
= e^{\delta}
{\displaystyle\prod_{\alpha\in R_+}}
\left( 1 - e^{-\alpha} \right) \in {\bf C}[H],
\]
where  
$\delta=\frac12 {\displaystyle\sum_{\alpha \in R_+}} \alpha \in P$.
\index{delta@$\delta$}
\begin{cor}\label{cor:skew}
Skew functions in ${\bf C}[H]$ are divisible by $\Delta$.
If we denote the set of $W$-skew Laurent polynomials by 
\index{che@${\bf C}[H]^{-W}$} ${\bf C}[H]^{-W}$,
then 
${\bf C}[H]^{-W} = \Delta{\bf C}[H]^W$.
\index{chw@$ {\bf C}[H]^W$}
\end{cor}
\pf 
Let $p \in {\bf C}[H]^{-W}$.
The previous proposition says that
$p \in (1-e^{-\alpha}) {\bf C}[H]$.
Since the algebra ${\bf C}[H]$ has the unique factorization property,
and $(1-e^{-\alpha})$ are coprime,
$p$ can be divided by $\Delta$.
\qed 
\begin{cor}
We put $\varepsilon(w) = \det_{\mathfrak a} w$.  
Then we have
\[
\Delta = \sum_{w \in W} \varepsilon(w) e^{w\delta}.
\]
\end{cor}
\pf  
Since the right hand side is skew, we have
\[\frac{1}{\Delta} \sum_{w \in W} \varepsilon(w) e^{w\delta}
\in {\bf C}[H]^W. \]
Moreover the leading term in the dominance ordering must be $1$.
\qed 

\index{k@$k=(k_\alpha)_{\alpha\in R}$}
Let $k_\alpha \in {\bf C}$ be $W$-invariant root labels,
that is, $k_\alpha = k_\beta$ 
if $\alpha,\beta$ are in the same $W$-orbit. 
We call $k=(k_\alpha)_{\alpha\in R}$ a multiplicity function on $R$.  
In this lecture we mainly consider real multiplicity functions and 
often assume that 
 $k_\alpha\geq 0$ for any $\alpha\in R$. 
\index{rho@$\rho(k)$}
We set
\[\rho(k)
=\frac12 {\displaystyle\sum_{\alpha \in R_+}}
 k_\alpha \alpha \in \mathfrak{h}^*.
\]
\index{Dunkl-Cherednik operator}
The hero of our story is the {\it Dunkl-Cherednik operator},
given by the following formula:
\begin{dfn}[Dunkl-Cherednik operator]\label{dfn:duop}
\index{txi@$T_\xi(k)$}
For $\xi \in \h$ define
\[
T_\xi(k) = \partial_\xi + 
\sum_{\alpha \in R_+} k_\alpha \alpha(\xi)
\frac{1}{1-e^{-\alpha}}(1-r_\alpha)
- \rho(k)(\xi).
\]
Here $\partial_\xi$ denote the  invariant
vector field on the torus $H$ corresponding to $\xi\in \h$.
\end{dfn}

\begin{remark}
By Proposition~1.1,
$T_\xi(k)$ maps ${\bf C}[H]$ to itself.
We may also think of $T_\xi(k)$ as an operator acting on 
other function spaces on $\mathfrak h$,
for example, holomorphic functions, 
or $C^\infty(A)$, 
or $C_c^\infty(A)$.
\end{remark}
\subsection{Commutativity} 
\begin{theorem} For any $\xi, \eta \in \mathfrak h$, we have
\[
[ T_\xi(k), T_\eta(k) ] = 0.
\]
\end{theorem}
\pf 
There are basically three proofs. 
A direct computation as in Dunkl's original paper, 
Cherednik's approach from 
conformal field theory (KZ equation),  
and Heckman's proof using orthogonality.
We give Heckman's proof here.

We introduce two important structures on ${\bf C}[H]$. 
In the rest of this section 
we assume $k_\alpha\geq 0$ for any $\alpha\in R$.  
First, we define the hermitian inner product
\index{(@$(f,g)_k$}
\[
(f,g)_k = 
{\displaystyle\int_T} f\, \bar{g}\, \delta_k\, dt,
\]
where the weight function is given by 
\index{deltak@$\delta_k$}
\[\delta_k
 = 
{\displaystyle\prod_{\alpha\in R_+}}
\left| e^{\alpha/2} - e^{-\alpha/2} \right|^{2k_\alpha}
=
{\displaystyle\prod_{\alpha\in R}}
\left| 1 - e^{\alpha} \right|^{k_\alpha}.
\]

Second, 
\index{<<@$\triangleleft$}
we introduce a partial ordering $\triangleleft$ on $P$ as follows :
$\lambda \triangleleft \mu$ if \index{<@$<$}
either $\lambda_+ < \mu_+$ in dominance ordering 
(with \index{lambdap@$\lambda_+$} $\lambda_+$ the 
unique dominant weight in $W \lambda$), or if 
$\lambda_+ = \mu_+$ and $\lambda > \mu$. 
This the last inequality is not 
a typographical error!
The following lemma explains the importance of the ordering and the 
inner product defined above:
\begin{lemma}
The operator $T_\xi(k)$ is upper triangular
with respect to $\triangleleft$, and $T_\xi(k)$ is 
symmetric with respect to $(\cdot, \cdot)_k$
if $\xi \in \mathfrak a$.
\end{lemma}
\pf
Using Proposition~\ref{prop-1.1},
we check that $T_\xi(k)$ is upper triangular
with respect to $\triangleleft$. The symmetry 
property is a simple direct computation left to the reader.
\qed

\begin{definition}
Define a basis
\index{ela@$E(\lambda,k)$}
$\{E(\lambda,k)  \,;\, \lambda \in P \}$ of ${\bf C}[H]$
by the following conditions.
\begin{itemize}
\item[(a)] $E(\lambda,k) 
= e^\lambda 
+ \displaystyle\sum_{\mu \triangleleft \lambda} c_{\lambda,\mu} e^{\mu}$.
\item[(b)] For any $\mu \triangleleft \lambda$,
$( E(\lambda,k), e^\mu )_k = 0$.
\end{itemize}
\end{definition}

Now we come back to the proof of Theorem~2.6.
$T_\xi(k) E(\lambda, k)$ also satisfies (a) and (b),
except that its expansion in (a) has leading term 
$(\tilde{\lambda}(\xi))e^{\lambda}$ for some $\tilde{\lambda}$.
The uniqueness shows that
\begin{equation}\label{eqn:ep}
T_\xi(k) E(\lambda, k) 
= \tilde{\lambda}(\xi) E(\lambda,k).
\end{equation}
Therefore $\{E(\lambda,k)\,;\,\lambda\in P\}$ 
diagonalize  simultaneously the 
Dunkl-Cherednik operators $T_\xi(k)$, 
hence these operators must mutually commute. \qed 

The eigenvalue $\tilde{\lambda}(\xi)$ can be calculated easily by 
Proposition~1.1:
\begin{cor}\label{cor:ev}
Define $\varepsilon\,:\,\R\rightarrow\{\pm 1\}$ by 
\[
\varepsilon(x) = \left\{
\begin{array}{ll}
1 & x>0 \\  -1 & x\le 0
\end{array}
\right.
\]
Given $\lambda \in P$, the eigenvalue in 
equation \ref{eqn:ep} 
is given by
\index{lamt@$\tilde\lambda$} 
\[
\tilde\lambda
= \lambda + \frac12 \sum_{\alpha \in R_+} k_\alpha
\varepsilon(\lambda(\alpha^\vee)) \alpha
= \lambda + w_\lambda^*(\rho(k)),
\]
where  \index{wls@$w_\lambda^*$} $w_\lambda^*$ is the longest 
element in $W$ sending
$\lambda_+$ to $\lambda$.
\end{cor}
\pf 
By Proposition~1.1, the eigenvalue $\tilde\lambda$ is given by
\begin{eqnarray*}
\tilde\lambda &=&
\lambda - \rho(k) +
  \sum_{\alpha \in R_+, \lambda(\alpha^\vee)>0} k_\alpha \alpha \\
&=& \lambda + \frac12 \sum_{\alpha \in R_+} k_\alpha
\varepsilon(\lambda(\alpha^\vee)) \alpha \\
&=& w_\lambda^* (\lambda_+ + \rho(k)) \\
&=& \lambda + w_\lambda^*(\rho(k)).
\end{eqnarray*}         
\qed 

Notice that 
the function $\varepsilon$ is not skew symmetric at $x=0$.
We can decompose $\mathfrak{a}^*$ in a non symmetric way in 
the disjoint ``chambers'' $C_w=\{\lambda\in \mathfrak{a}^*\mid
\lambda(\alpha^\vee)>0\ \forall \alpha\in R_+\cap w(R_+)\ {\rm{ and }}
\ \lambda(\alpha^\vee)\leq 0\ \forall \alpha\in R_+\cap w(R_-)\}
$ 
(with $w$ traversing $W$) 
which lie between $w(\mathfrak{a}^*_+)$ 
and $\overline{w(\mathfrak{a}^*_+})$. The map $\lambda\to\tilde\lambda$
restricted to $C_w$ is a translation by the vector $w(\rho(k))$. 
So the chambers $C_w$ are shifted apart from each other by this 
map, and the joint spectrum of the $T_\xi(k)$ operators on 
${\bf C}[H]$ is obtained by applying this map to the lattice $P$.
\begin{cor}
$\{ E(\lambda,k) \,;\, \lambda \in P \}$ is an orthogonal basis
of ${\bf C}[H]$
{\rm (}assuming still that $k_\alpha \ge 0$ for any $\alpha${\rm )}.
\end{cor}
\begin{pf} The eigenvalues $\tilde{\lambda}$ are mutually distinct.
\qed
\end{pf}


``Macdonald theory'' is concerned with these polynomials
$E(\lambda,k)$ and their further properties,
for example,
the computation of their $L^2$ norm with respect to
$(\cdot, \cdot)_k$,
and their value at $e \in H$.
To attack these problems effectively,
we must investigate the algebraic structures 
attached to the $T_\xi(k)$.
This is the main subject of the next three sections.

\newpage
%
\section{Degenerate double affine Hecke algebra}

The results in this section are due to Ivan Cherednik, see \cite{C1}, \cite{C2}.

\index{affine Weyl group}
The {\it affine Weyl group} \index{wa@$W^a$} $W^a$ is the group 
acting on $\mathfrak h^\ast$, generated by 
the reflections \index{ra@$r_a$} $r_a$, $a = [\alpha ^\vee,n] \in  
R^\vee+\mathbf{Z}\subset S(\mathfrak{a})$, defined by 
\[
r_a(\lambda ) = r_{[\alpha ^\vee,n]} (\lambda ) 
= \lambda - (\lambda (\alpha ^\vee ) + n) \alpha .
\]
We shall often write $a=\alpha^\vee +n$ as an 
element of $S(\mathfrak{a})$ instead of $[\alpha^\vee, n]$. 
 In particular, this group contains all translations in $Q$, 
since for any $\alpha \in R$, 
\[
  r_{\alpha ^\vee}r_{[\alpha^\vee,1]} = 
r_{[-\alpha^\vee,1]}r_{\alpha^\vee} = t_\alpha.
\] 
\index{ta@$t_\alpha$}
In fact, one has $W^a=W \ltimes Q$, 
the semidirect product of $Q$ by $W$. 
This is a Coxeter group of affine type, 
if we take the set of simple reflections for $W^a$ equal to  
$\{r_0, r_1, \cdots ,r_n\}$, with $r_i = r_{a _i}$, 
\index{a0@$a_0$} $a_0 = [-\theta^\vee,1]$, 
and $a_i = \alpha_i^\vee$, $i>0$. 
Here \index{theta@$\theta$} $\theta$ denotes the 
unique {\it highest short root}. 


The affine positive roots are \index{rap@$R^a_+$}$R^a_+  =
R_+ \cup (R+\mathbf{Z}_{>0})$, 
and the corresponding set of simple roots is denoted 
\[S^a=\{a_0,a_1,\cdots,a_n\}.\] 
The fundamental alcove \index{c@$C$} $C$ is 
\[C=\{ \lambda \in \mathfrak a^\ast \,;\, 
\lambda(a_i)>0, \  \ i=0,1,\cdots,n \}.
\] 
Then $\overline{C}$ is a fundamental domain for the action of $W^a$.

\index{extended affine Weyl group}
We shall work with  \index{we@$W^e$} $W^e = W \ltimes P$, 
the {\it extended affine Weyl group}. 
This is not a Coxeter group in general, 
but $W^a \triangleleft W^e$ and if 
\index{omega@$\Omega$} 
\[\Omega
= \{ \omega \in W^e \,;\,  \omega(C)=C   \},
\] 
then $\Omega \cong P/Q$, and 
\[ 
W^e =  W^a\ltimes \Omega .
\]
Clearly $\omega \in \Omega$ 
defines a permutation of the set $S^a$.


By {\it duality} the action of $W^e$ on 
$\mathfrak h^\ast$ via affine transformations gives 
rise to a representation of $W^e$ on the symmetric 
algebra $S(\mathfrak h)$ of $\mathfrak h$ 
(viewed as polynomial functions on $\mathfrak h^\ast$). 
Notice that \index{snh@$S_n(\mathfrak h)$} $S_n(\mathfrak h)$ 
(the part of $S(\mathfrak h)$ of degree $\leq n $) 
is stable under this action; for  $n=1$  
this action gives the {\it reflection representation} 
of $W^e$ on $\mathfrak h \oplus 
\mathbf{C}$, explicitly given by:
\[
 r_{[\alpha^\vee,n]} [\xi,u] = [\xi,u]-\alpha(\xi)[\alpha^\vee,n],
\]
and
\[
t_\lambda [\xi,u] = [\xi,u-\lambda(\xi)],
\]
where $[\xi,u] (\lambda) = \lambda(\xi)+u$. 
If $ p \in S(\mathfrak h)$, and $w \in W^e$, 
then write \index{pw@$p^w$} $p^w(\lambda)= p(w^{-1}\lambda)$.


Since we need to understand precisely the relation 
$\Omega \cong P/Q$ we introduce the following notion.

\begin{df}
\index{minuscule weight}
An element in $\overline{C}\cap P\setminus\{0\}$ is called 
a {\it minuscule weight}.
\end{df}

\begin{prop}
Let $\{\lambda_1,\lambda_2,\cdots,\lambda_n \}$ denote 
 the set of fundamental weights for the simple system 
$\{r_1,r_2,\cdots,r_n \}$ 
and  $\theta^\vee = \sum_{i=1}^n n_i a_i$ the maximal coroot. 
Put $O^\ast = \{i \in \{1,2,\cdots,n\} \,;\, \ n_i=1 \}$. 
Then $\overline{C}\cap P\setminus\{0\} 
= \{ \lambda_i\,;\, \ i\in O^\ast \}$.
\end{prop}
\pf 
\index{lambdai@$\lambda_i$}
Obviously $\overline{C}\cap P\setminus\{0\} \supset \{ \lambda_i
\,;\,  i\in O^\ast \}$. 
In the other direction we argue as follows. If $\lambda\in\overline{C}
\cap P\setminus\{0\}$ then $\lambda(\theta^\vee)=1$.
Write $\lambda=\sum_{i=1}^n m_i\lambda_i$, and notice that $m_i \in 
\mathbf{Z}_{\geq 0}$ and that $\lambda_i(\theta^\vee)\in\Z_{>0}$. Hence from 
\[
\lambda(\theta^\vee) = \sum_{i=1}^n m_i \lambda_i(\theta^\vee)=1,
\]
it follows that there exists an $i$ such that 
$m_i = \lambda_i(\theta^\vee) =1$ and $m_j =0$ ( for $i \neq j$). 
Thus $\lambda=\lambda_i$ and $i\in O^\ast$.
\qed 

For $r \in O^\ast$, 
let \index{omegar@$\omega_r$} 
$\omega_r= t_{\lambda_r}w_{\lambda_r}w_0 \in W^e$,  
where $w_{\lambda_r}$ is the longest element 
in the parabolic subgroup $W_{\lambda_r}$ of $W$ 
generated by $\{r_1,\cdots,r_{r-1},r_{r+1},\cdots,r_n \}$ 
(the stabilizer of $\lambda_r$)  
and $w_0$ is the longest element in $W$. The 
parabolic subsystem of roots that corresponds to 
$W_{\lambda_i}$ is denoted by 
$R_{\lambda_i}$. Its basis of simple roots is 
$\{\alpha_1,\dots,\alpha_{r-1},\alpha_{r+1},\dots,
\alpha_n\}$.


\begin{prop}
$\Omega= \{\omega_r \in W^e \,;\, r\in O^\ast \} \cup \{id_\mathfrak{a^*} \}$. 
In particular the set of all minuscule weights is a complete 
set of representatives of $P/Q \setminus \{0\}$. 
\end{prop}
\pf 
Let $\omega \in W^e$ such that $\omega(C)=C$. 
Then $\omega(S^a)=S^a$, where $S^a=\{a_0=1-\theta^\vee,a_1,\cdots,a_n \}$. 
If $\omega(a_0)=a_0$, 
then $\omega(\{\alpha_1,\cdots,\alpha_n \})=\{\alpha_1,\cdots,\alpha_n \}$, 
therefore $\omega=id_\mathfrak{a^*}$ by simple transitivity of the 
action on chambers of $W$. Hence we may and will label $\omega\in\Omega$ uniquely 
by the index $r\in\{0,\dots,n\}$ such that $\omega_r(a_0)=a_r$. 
Now let $r\in \{1,2,\cdots,n\}$, and write  
$\omega_r=t_{\mu_r}w_r$. Then $w_r(\theta^\vee)=-a_r$ 
and $\mu_r \in \overline{C}\cap P\setminus\{0\}$. 
Hence $\mu_r$ is a 
minuscule fundamental weight and $\mu_r(a_r)=1$. In other words, 
it is the fundamental 
weight $\lambda_r$ of $a_r$. 
Because $\omega_r^{-1}= \  w_r^{-1}t_{-\lambda_r}$ we have   
$w_0 w_r^{-1}(\lambda_r) \in \overline{C}$. Hence  
$w_0 w_r^{-1}(\lambda_r)=\lambda_r$.
Moreover, for $i \neq r$ we have 
$w_0w_r^{-1}(\alpha_i) = w_0(\alpha_j) \in R_-$  
for some $j \in \{1,2,\cdots,n\}$.
Therefore we have $w_0w_r^{-1} = w_{\lambda_r}$, hence 
$w_r = w_{\lambda_r}w_0$.
 
Vice versa, let $\lambda_r$ be a minuscule fundamental weight. 
Since $w_0\mu \in -C$ for $\mu \in C$ and $w_{\lambda_r}(a_i) 
\in R_{\lambda_r,-}$ $(i\neq 0,r)$, we have 
\[
\omega_r \mu(a_i) \ = \lambda_r(a_i)+\  w_0 
\mu(w_{\lambda_r}(a_i))= w_0 \mu(w_{\lambda_r}(a_i))  \ >0. 
\]   
Since $\theta^\vee \geq w_{\lambda_r}(a_r)$ and $w_0\mu (\theta^\vee)>-1$, 
we have
\[
\omega_r \mu (a_r) =\lambda_r(\alpha_r)+\  w_0 \mu(w_{\lambda_r}(a_r)) 
=  1+ w_0\mu(w_{\lambda_r}a_r)>1+w_0\mu(\theta^\vee)>0.
\] 
On the other hand, $w_{\lambda_r}(\theta^\vee) \in R_+^\vee$ 
and $\lambda_r$ is a minuscule weight, thus  
\[
\omega_r\mu(\theta^\vee) = 1+w_0\mu (w_{\lambda_r}\theta^\vee) <1.
\]
Thus we have $\omega_r C \subset C$, that is  $\omega_r \in \Omega$.
The map $O^\ast \ni r \to \omega_r \in \Omega$ is injective 
since $\omega_r(0)=\lambda_r$.   
\qed 

\begin{cor} (of proof)
If $\lambda_r$ is a minuscule weight, then $\omega_r(1-\theta^\vee)= a_r$.
\end{cor}


\begin{df}{\rm (Cherednik)}\label{def:hecke}
\index{degenerated extended double affine Hecke algebra}
The {\it degenerated extended double affine Hecke algebra}  
\index{her@$\mathbf{H}^e(R_+,k)$} $\mathbf{H}^e(R_+,k)$ is 
the unique associative algebra 
over $\mathbf{C}$ such that \\
\smallskip
$(1)$ $\mathbf{H}^e (R_+,k) 
\cong S(\mathfrak h) \otimes \mathbf{C}[W^e]$ as 
vector space over $\mathbf{C}$, \\
\smallskip
$(2)$ $S(\mathfrak h) \ni p \mapsto p\otimes e 
\in \mathbf{H}^e(R_+,k)$, and $\mathbf{C}[W^e] \ni w 
\mapsto 1\otimes w \in  \mathbf{H}^e(R_+,k)$ are algebra homomorphisms, \\
\smallskip
$(3)$ $(p\otimes e)(1\otimes w) = p\otimes w$. \\
Write $p\cdot w$, or $p w$ instead of $p\otimes w$ from now on. \\
\smallskip
$(4)$ $r_i\cdot p-p^{r_i}\cdot r_i =-k_i (p-p^{r_i})/a_i$, 
$(i=0,1,\cdots,n),$ 
where $k_0 =k_\theta$. \\
\smallskip
$(5)$ $\omega \cdot p = p^\omega \cdot \omega$ for all 
$\omega \in \Omega$.     

\end{df}


\begin{thm}{\rm (Cherednik)}
Let \index{a@$\mathbf{A}$} $\mathbf{A}$ denote a subalgebra of 
$\mathrm{End}(\mathbf{C}[H])$ generated by $e^\lambda$ $(\lambda \in P)$, 
$w \in W$, and $T_\xi(k)$ $(\xi \in \mathfrak{h} )$. 
Then 
\index{pi@$\pi$}
\[\pi
:
W^e \ni t_\lambda w \mapsto e^\lambda w \in  \mathrm{End}(\mathbf{C}[H])
\]   
and
\[
\pi:\mathfrak{h} \ni \xi \mapsto T_\xi(k) \in  \mathrm{End}(\mathbf{C}[H])
\]
extend to a representation of $\mathbf{H}^e(R_+,k)$ on 
$\mathbf{C}[H]$, and $\mathbf{H}^e(R_+,k)$ is isomorphic 
to $\mathbf{A}$ via $\pi$.
\end{thm}
\pf 
We need to check (4) and (5), the other points being obvious. 

First notice that $\pi:W^e \to \mathrm{End}(\mathbf{C}[H])$, 
and $\pi:S(\mathfrak{h}) \to \mathrm{End}(\mathbf{C}[H])$ 
are well defined. We can check by simple direct 
computation that  $T_\xi(k)$ 
and $r_i$ $(i=1,2,\cdots,n)$ satisfy the relation $(4)$. 
The case $r=0$ requires a bit of special care: 
put
\index{sxi@$S_\xi(k)$} 
\[
S_\xi(k)= \partial_\xi + 
\frac{1}{2}\sum_{\alpha \in R_+} k_\alpha \alpha(\xi)
\frac{1+e^{-\alpha}}{1-e^{-\alpha}}(1-r_\alpha),
\]
\index{Dunkl-Heckman operator}
This operator is called  the {\it Dunkl-Heckman operator}. Define 
$u_\xi(k)$ by $T_\xi(k) = S_\xi(k) - u_\xi(k)$, then 
\index{uxi@$u_\xi(k)$}
\[u_\xi(k)
 = \frac12 {\displaystyle\sum_{\alpha \in R_+}} 
k_\alpha \alpha(\xi) r_\alpha .
\]

The operator $S_\xi(k)$ is independent of the choice of a 
positive system $R_+$ of $R$ and $wS_\xi(k)w^{-1} = S_{w\xi}(k)$ 
for all $w\in W$, $\xi \in \mathfrak{h}$ 
(but $\{S_{\xi}\,;\,\xi\in\mathfrak{h}\}$ is not commutative). 
We leave it to the reader to verify by direct computation that 
\[
\pi(r_0) S_\xi(k)\pi(r_0) =
S_{r_0(\xi)}(k) - \frac{1}{2}\sum_{\alpha \in R_+}k_\alpha \alpha(r_0(\xi))
\left\{
\frac{(1-e^{\theta(\alpha^\vee)\alpha})(1+e^\alpha)}{1-e^\alpha} 
\right\} r_\alpha,
\]
and
\[
\pi(r_0) u_\xi(k) \pi(r_0) = -\frac{1}{2}\sum_{\alpha \in R_+} 
k_\alpha \varepsilon(\theta(\alpha^\vee))\alpha(r_0(\xi))e^{\theta(\alpha^\vee)\alpha}r_\alpha.
\]
Using that $\theta(\alpha^\vee)=0$ or $1$ we now check the desired relation 
$\pi(r_0) T_\xi(k)\pi(r_0) = 
T_{r_0(\xi)}(k) + k_0 \theta(\xi)\pi(r_0)$.


Let's look at relation (5). For the minuscule fundamental weight 
$\lambda_r$ of a simple root $\alpha_r$, 
we put \index{pir@$\pi_r$} $\pi_r :\mathbf{C}[H] \to \mathbf{C}[H]$, 
$\pi_r= \pi(\omega_r)=e^{\lambda_r}w_r$. Straightforward computations show:
\[
\pi_r S_\xi(k) \pi_r^{-1} = S_{\omega_r(\xi)}(k)+
\frac{1}{2}\sum_{\alpha \in R_+}k_\alpha \alpha(\xi)
(r_{w_r\alpha}-r_{\omega_r\alpha})(-\lambda_r(w_r\alpha^\vee))
\]%
and
\[
\pi_r u_\xi(k)\pi_r^{-1} = -\frac{1}{2}\sum_{\alpha\in R_+} k_\alpha \alpha (\xi) r_{\omega_r \alpha},
\]
hence
\begin{eqnarray*}
\pi_r T_\xi(k)\pi_r^{-1} &=& S_{\omega_r(\xi)}(k) + 
\frac{1}{2}\sum_{\alpha\in R_+}\varepsilon(-\lambda_r(w_r\alpha^\vee))
k_\alpha \alpha (\xi) r_{w_r \alpha}  \\
                         &=& T_{\omega_r(\xi)}(k).
\end{eqnarray*}


Finally we show that $\pi$ is an isomorphism. 
Obviously $\pi$ is surjective. 
Suppose that $\sum_{w\in W} p_w(T(k))w=0$ in $\mathbf{A}$. 
If we write $\sum_{w\in W} p_w(T(k))w = \sum_{w\in W}D_w w$, 
then $D_w=0$ for all $w \in W$. On the other hand, let $w'$ 
be such that the degree of $p_{w'}$ is maximal and let $q$ 
denote its highest degree part. Then the highest order part of 
$D_{w'}$ equals $\partial_p$, hence $q=0$. Consequently, 
$p_w=0$ for all $w \in W$.   
\qed

 
We can give a more intrinsic definition of the model representation: 

\begin{df}{\rm (\cite{Dri}, \cite{Lu})}
$\mathbf{H}(R_+,k) \cong S(\mathfrak{h})\otimes \mathbf{C}[W] 
\subset \mathbf{H}^e(R_+,k)$ is called the 
\index{degenerate affine Hecke algebra}
{\it degenerate affine Hecke algebra} 
\index{graded affine Hecke algebra}
or {\em graded affine Hecke algebra}.
\end{df}

\begin{df}
We can define a one dimensional representation of $\mathbf{H}(R_+,k)$ by 
\[
 \left\{
  \begin{array}{ll}
   \xi \cdot 1 = - \rho(k)(\xi)  1
   & (\xi \in \mathfrak{h}) \\
   w \cdot 1 = 1
   & (w\in W). 
  \end{array}
 \right.
\]
This representation is called the 
{\it trivial representation} of $\mathbf{H}(R_+,k)$, which 
we denote by \index{triv} triv.
\end{df}

\begin{thm}
The representation $\pi$ is isomorphic to the induced 
representation 
$\mathrm{Ind}_{\mathbf{H}(R_+,k)}^{\mathbf{H}^e(R_+,k)}(\mbox{\rm triv})$.
\end{thm}
\pf 
For $1 \in \mathbf{C}[H]$, 
$T_\xi(k)\cdot 1 = -\rho(k)$ and $w\cdot 1=1$. 
Hence there exist a unique epimorphism $\varphi 
:\mathrm{Ind}_{\mathbf{H}(R_+,k)}^{\mathbf{H}^e(R_+,k)}(\mbox{triv}) 
\to \pi$ such that $\varphi(1)=1$. On the other hand,
 as a $\mathbf{C}[H]$ module, 
$\mathrm{Ind}_{\mathbf{H}(R_+,k)}^{\mathbf{H}^e(R_+,k)}(\mbox{triv})$ 
is isomorphic to the left 
regular representation of $\mathbf{C}[H]$. 
Hence, as a $\mathbf{C}[H]$ module, 
$\mathrm{Ind}_{\mathbf{H}(R_+,k)}^{\mathbf{H}^e(R_+,k)}(\mbox{triv}) 
\cong \pi$. 
Therefore, as a $\mathbf{H}^e(R_+,k)$ module, 
$\mathrm{Ind}_{\mathbf{H}(R_+,k)}^{\mathbf{H}^e(R_+,k)}(\mbox{triv})$ is 
isomorphic to $\pi$ via $\varphi$.
\qed 
\newpage     
%
\section{Intertwiners}
The intertwining operators between minimal principal series representations 
of (graded) affine Hecke algebras are built from certain intertwining elements 
of these algebras. This is a main topic of study in the representation theory 
of Hecke algebras. In this section we will extend this construction to the
double affine situation, and discuss the basic applications to Macdonald
theory. The ideas in this section are mainly due to Ivan Cherednik.

\subsection{Intertwining elements in the degenerate double affine Hecke
algebra} 

In the degenerate graded Hecke algebra there exist elements $I_w$
for $w\in W^e$  with the property that the conjugate inside
$\mathbf{H}^e(R_+,k)$ of an  element $p\in S(\mathfrak h)$ by $I_w$ is equal
to $p^w$.  These elements are called 
``intertwiners'', because they give rise to intertwining maps between minimal 
principal series modules.  In our context this means that we find operators $\pi(I_w)$ 
which map solutions of \ref{eqn:ep} to solutions of \ref{eqn:ep} with spectral 
parameter $w\lambda$.


\begin{df}
\index{ii@$I_i$} 
\[
I_i
= r_i a_i + k_i \  \in \mathbf{H}^e(R_+,k) \  \  (i=0,1,\cdots,n)
\]
\end{df}

\begin{thm}
$($a$)$ $I_i^2 = k_i^2 - a_i^2 $. \\
$($b$)$  $I_i  p = p^{r_i} I_i $ $\forall p\in S(\mathfrak h)$.\\
$($c$)$  $I_i I_j I_i \cdots = I_j I_i I_j \cdots $  \\
with $m_{ij}$ factors on both sides. Here $m_{ij}$ denotes 
the order of the element $r_i r_j \in W^a$.
 \\
$($d$)$  Assume that $k_\alpha\geq 0$ for all $\alpha\in R$. 
Then we have $(I_i f,g)_k = -(f,I_i g)_k $  for all $i = 0,1,\cdots,n$.
\end{thm}
\pf 
$($a$)$ and $($b$)$ are trivial reformulation of (4) in 
Definition~\ref{def:hecke}, 
and $($d$)$ follows directly from the symmetry of $T_\xi(k)$. 
Statement (c) is equivalent with  the following; 
if we have two reduced expressions 
$r_{i_1} r_{i_2} \cdots r_{i_n} = r_{i'_1} r_{i'_2} \cdots r_{i'_n}$ 
for $w$, 
then $I_{i_1}I_{i_2} \cdots I_{i_n} = I_{i'_1}I_{i'_2}\cdots I_{i'_n}$.
For a reduced expression  
$r_{i_1} r_{i_2} \cdots r_{i_n}$, 
we put \index{iw@$I_w$} $I_w =I_{i_1}I_{i_2} \cdots I_{i_n}$. 
Notice that we can write 
\[
I_w = w  \prod_{a \in R_+^a, w(a) \in R_-^a} a  + \sum_{w' < w} p_{w,w'} w',
\]
where $p_{w,w'} \in S(\mathfrak h)$, thus, 
if we allow rational coefficients, we also have 
\[ 
I_w = w  
\prod_{a \in R_+^a, w(a) \in R_-^a} a  + \sum_{w' < w} r_{w,w'} I_{w'}.
\]
The top coefficient is independent of the reduced expression for $w$; 
so if $I_w$ and $I'_w$ are different, then the difference 
$I''_w = I_w -I'_w $ is of the form $\sum_{w'<w} r'_{w,w'}I_{w'}$ 
and also have intertwining property $I''_w p = p^w I''_w$ 
$(p \in S(\mathfrak h))$. Thus we have  $I''_w = 0$.
\qed


By the above theorem, we can define $I_w$ for $w \in W^a$ as follows; if $w=
r_{i_1} r_{i_2} \cdots r_{i_n}$ is a reduced expression for $w$, then we put 
\[ 
I_w = I_{i_1}I_{i_2} \cdots I_{i_n} .
\] 

Obviously, we also have $\omega I_i = I_j \omega $ if 
$\omega \in \Omega$ and $\omega r_i = r_j \omega $. 
Hence we may also use $\Omega$ to build intertwiners for 
arbitrary elements of $W^e$:

\begin{df}
For a reduced expression 
$w = \omega r_{i_1} r_{i_2} \cdots r_{i_n}$ for $w \in W^e$, 
we define the {\it general intertwiner} 
$I_w \in \mathbf{H}^e(R_+,k)$ for $w$ by
\[
I_w = \omega I_{i_1}I_{i_2} \cdots I_{i_n} .
\]
\end{df}

\begin{cor}\label{cor:rod}
For $w \in W^e$ we have 
\[ 
I_w(1) = d(w,k) E(w(0),k), \] 
where
\index{dwk@$d(w,k)$}
\[
d(w,k)
= \prod_{a \in R_+^a \cap w^{-1} R_-^a} a(-\rho(k)).
\]
\end{cor}    

\begin{remark}\label{rem:coc}
The equality $I_w = \omega I_{i_1}I_{i_2} \cdots
I_{i_n}$ is true only if the expression $w=\omega r_{i_1} r_{i_2} \cdots
r_{i_n}$ is reduced. Denote by $I_w(\lambda)$ the right evaluation of $I_w$ at
$\lambda$. In other  words, $I_w(\lambda)$ is the element of $\C[W^e]$ defined
by  
\[
I_w(\lambda) = \omega I_{i_1}(r_{i_2}\dots r_{i_n}\lambda)I_{i_2}( r_{i_3}\dots
r_{i_n}\lambda)\cdots I_{i_n}(\lambda)
\]
with $I_{i}(\lambda)=\lambda(a_i)r_i+k_i$.
If we normalize these elements of $\C[W]$ as follows:
\[
{\tilde I}_w(\lambda)=\frac{I_w(\lambda)}{\prod_{\alpha\in R^a_+\cap
w^{-1}(R^a_-)}(\lambda(a)+k_a)}
\]
then the ${\tilde I}_w(\lambda)$ behave as a $W^e$ cocycle:
\[
{\tilde I}_{ww^\prime}(\lambda)={\tilde I}_w(w^\prime\lambda)
{\tilde I}_{w^\prime}(\lambda)
\]
for all $w$ and $w^\prime$.
\end{remark}

\subsection{Application: Macdonald's conjectures}

Intertwiners can be used to verify the Macdonald's norm and evaluation 
conjectures. This was not the first proof of these conjectures, but it
is the most natural proof at this point. (The original approach was based on
the so-called shift principle, which will be discussed in the next section.)

\begin{df}
For $w \in W$ we put
\[\delta_w(\alpha)= 
  \left\{
   \begin{array}{ll}
    0 & 
    \mbox{ if $\alpha \in w^{-1}R_+$} \\
    1 &
    \mbox{ if $\alpha \in w^{-1}R_-$} 
   \end{array}
  \right. .
\]
We define meromorphic functions \index{csw@$c^\ast_w(\lambda,k)$} 
$c^\ast_w$ and 
\index{ctw@$\tilde{c}_w(\lambda,k)$} $\tilde{c}_w$ in $\lambda$, $k$ by 
\begin{equation}\label{eqn:cast}
c^*_w(\lambda,k) = 
\prod_{\alpha\in R_+}
\frac{\Gamma(-\lambda(\alpha^\vee)-k_\alpha+\delta_w(\alpha))}
{\Gamma(-\lambda(\alpha^\vee)+\delta_w(\alpha))},
\end{equation}
\begin{equation}\label{eqn:ctilde}
\tilde{c}_w(\lambda,k) =
\prod_{\alpha\in R_+}
\frac{\Gamma(\lambda(\alpha^\vee)+\delta_w(\alpha))}
{\Gamma(\lambda(\alpha^\vee)+k_\alpha+\delta_w(\alpha))}.
\end{equation}
In particular we put \index{ct@$\tilde{c}(\lambda,k)$} $\tilde{c}=\tilde{c}_e$.
\end{df}

For $\lambda \in P_+$, 
we put  \index{wl1@$W_\lambda$} 
$W_\lambda=\{w\,;\,w\lambda=\lambda\}$ and 
\index{wl2@$W^\lambda$}
$W^\lambda=
\{w\,;\,l(ww')\geq l(w) \mbox{ for all }w'\in W_\lambda\}$. 
Let $w_\lambda$ denote the longest element in 
$W_\lambda $. 
\begin{thm}\label{thm:norm}
Assume that $k_\alpha\geq 0$ for all $\alpha\in R$. 
For $\lambda \in P_+$
 and  $w \in W^\lambda$, we have 
\[
\Vert E(w\lambda,k) \Vert_k^2 = 
\frac{c^\ast_{w {w_\lambda}}(-(\lambda+\rho(k)),k)}
{\tilde{c}_{w{w_\lambda}}(\lambda+\rho(k),k)}, 
\]
and 
\[
E(w\lambda,k)(e) = 
\frac{\tilde{c}_{w_0}(\rho(k),k)}
{\tilde{c}_{w{w_\lambda}}(\lambda+\rho(k),k)} .
\]
\end{thm}
\pf
Use Corollary \ref{cor:rod}.
\qed

%
%
\subsection{Jack Polynomials}
When $R$ is of type $A_n$, Knop and Sahi used this approach to verify 
the integrality and positivity conjecture for Jack polynomials 
(also in the nonsymmetric case). 

\begin{thm}{\rm (F.Knop and S.Sahi \cite{KS})}
For a partition $\lambda$  of $n$ let $m_i(\lambda)$ be the number 
of parts which are equal to $i$ and let 
$u_\lambda =\prod_{i\geq1} m_i(\lambda) ! $. 
If the Jack polynomial \index{jlam@$J_\lambda (x;\alpha)$} 
$J_\lambda (x;\alpha)$ has a expansion 
\[
J_\lambda (x;\alpha)= \sum_{\nu \geq 0} v_{\lambda,\nu}(\alpha) m_\nu (x)
\]
by monomial symmetric functions 
$m_\nu$ $(\nu$\,:\,partition of $n)$, 
then all functions 
$\tilde{v}_{\lambda, \nu}= u_\lambda^{-1} v_{\lambda ,\nu}(\alpha)$ 
are polynomials in $\alpha$ with positive integral coefficients.
\end{thm}

Here, in terms of our notations, $\alpha$ is the 
inverse of the multiplicity $k$ and
\[
J_\lambda (x;\alpha) = \prod_{b \in \lambda } c_\lambda (b) 
\frac{1}{|W_\lambda|} \sum_{w\in W}  E^w(\lambda,k),
\]
where, for $\lambda$ and $b=(i,j) \in \lambda$; 
a box in $\lambda$, $c_\lambda(b)=\alpha(\lambda_i -j)+(
\mbox{leg}(b)+1)$.

\begin{remark}         
In fact Knop and Sahi proved a stronger result, namely   
a combinatorial formula for the Jack polynomial. 
\end{remark}

\newpage
%
\section{The shift principle}
In the previous section we introduced operators that act 
on the spectral parameter $\lambda$ of \ref{eqn:ep}. In this section 
we will study operations on the multiplicity parameter $k$. There exist 
so-called shift operators that induce translations in a certain lattice 
in the parameter space $K$. The most fundamental example of this kind 
of operator is already sufficient to prove Macdonald's constant term and
evaluation  conjectures, and therefore we will restrict ourselves to the
discussion of this simplest example of a shift operator. 

It is remarkable that these shift operators act naturally on 
the $W$ symmetrizations of solutions of \ref{eqn:ep}, rather than on the
solutions themselves. However, on the solution space of \ref{eqn:ep}, 
symmetrization for the action of $W$ is invertible by a differential operator.
This will become clear in the section on the KZ equation (see
Remark \ref{rem:inv}). 

The $W$ symmetrizations of solutions of 
\ref{eqn:ep} are eigenfunctions of an important system of commuting 
differential operators that will play the leading part in the next section.
This system is called the hypergeometric system of differential equations. In
the section  on the KZ equations we shall see that this system is generically
equivalent to \ref{eqn:ep} (Matsuo's theorem), but it represents a different
point of view (somewhat like spherical representations versus principal series
representations).

When considering these hypergeometric differential operators, yet 
another symmetry in the parameter space $K$ arises naturally. This is the 
reflection symmetry $k_\alpha^\prime=1-k_\alpha$, and this will also be
discussed in this section. 

\subsection{Translation symmetry in the multiplicity parameter}

In this section we use the notation 
\index{h@$\mathbf{H}$} $\mathbf{H}  = \mathbf{H}(R_+,k)$ for the degenerate
affine  Hecke algebra.  Here $k$ is a multiplicity  such that $k_\alpha \geq
0$ for all $\alpha \in R$. 

\begin{lm}\label{lm:cent}
\index{zh@$Z(\mathbf{H})$} $Z(\mathbf{H}) = S(\mathfrak h)^W $.
\end{lm}
\pf 
The following formula  
can be checked by induction on the length of $w$: 
\begin{equation}\label{eqn:wxi}
w\cdot \xi \cdot w^{-1} = w(\xi) + 
\sum_{\alpha \in R_+\cap wR_-} k_\alpha \alpha(w\xi)r_\alpha. 
\end{equation}
From this formula one deduces easily that $Z(\mathbf{H})\subset S(\mathfrak h)$.
Then one may use Definition~\ref{def:hecke} (4) to prove the result. 
 \qed 

\begin{df}
Let us define a subspace \index{mla@$M(\lambda,k)$} 
$M(\lambda,k)$ of $\mathbf{C}[H]$  by
\[
M(\lambda,k) = \{ f\in \mathbf{C}[H] \,;\,
p(T_\xi(k)) f = p(\lambda) f, p \in S(\mathfrak{h})^W \} .
\]
\end{df}

\begin{prop}  
For all $\lambda \in \mathfrak{h}^\ast$ we have;
\[
M(\lambda,k) =
  \left\{
   \begin{array}{ll}
   \mbox{\rm Span}\{ E(\nu,k) \}_{\nu \in W\bar{\lambda} }  & 
    \mbox{ if $\exists\bar{\lambda} \in P_+$ s.t. 
$\lambda \in W(\bar{\lambda}+\rho(k))$}, \\
    \{ 0 \} &
    \mbox{otherwise. } 
   \end{array}
  \right.
\]
\end{prop}
\pf
This follows from Corollary~\ref{cor:ev}. 
\qed 

\begin{cor}
$M(\lambda,k)$ is a module over $\mathbf{H}$
\end{cor}

As a module for $\mathbf{C}[W]$, $M(\lambda,k)$ 
is independent of $k$, of course, so $M(\lambda,k) 
\cong \mathbf{C}[W/W_\lambda] =\mathbf{C}[W^\lambda]$. 
In particular, there is a {\it unique} $W$-invariant element  
up to a scalar multiple.

\begin{df}
\index{Jacobi polynomial}
For $\lambda \in P_+$, 
the {\it Jacobi polynomial} \index{pl@$P(\lambda,k)$}
$P(\lambda,k)\in M(\lambda,k)$ 
is defined by 
\[
P(\lambda,k) = \sum_{w\in {W^\lambda}} E^w (\lambda,k),
\]
where $E^w$ denote the function on $T$ defined by $E^w(t)=E(w^{-1}t)$. 
Then it is of the form
\[P(\lambda,k)=
\sum_{\nu \in P_+, \nu \leq \lambda} c_{\lambda,\nu}(k) m_\nu ,                  \  \  c_{\lambda,\lambda}(k) = 1 .
\]
\end{df}

If $\lambda$ is {\it regular} in $P_+$, $M(\lambda,k)$ also contains a
one-dimensional skew-invariant subspace, and we can define a
skew-invariant function
\index{p-@$P^-(\lambda,k)$} 
\[P^-(\lambda,k)
= \sum_{w\in W} \varepsilon(w) E^w (\lambda,k).
\]
The next theorem is the heart of the ``shift principle''. 
It is a direct generalization of Weyl's character formula.
\begin{thm}{\rm (Generalized Weyl character formula)}
\[
P^-(\lambda+\delta,k) = \Delta P(\lambda,k+1)
\]
or
\[
P(\lambda,k+1) = 
\frac{P^-(\lambda+\delta,k)}{\Delta} = 
\frac{P^-(\lambda+\delta,k)}{P^-(\delta,k)}.
\]
\end{thm}
\pf
The assertion follows directly from the divisibility 
(Corollary \ref{cor:skew}) 
of skew polynomials by $\Delta$ and the definition of 
the $E(\lambda,k)$ using orthogonality. 
\qed 

It is not difficult to show that $M(\lambda,k)$ is 
{\it irreducible} as $\mathbf{H}$-module. 
Consequently, the shift principle is effective to understand  
properties of $M(\lambda,k)$ if $k_\alpha \in 
\mathbf{Z}_{>0}$ for all $\alpha \in R$, because it reduces 
everything to the trivial situation of $M(\lambda+\rho(k),0)$, via 
induction on $k$. 
For example we can prove Theorem \ref{thm:norm} in this way.

\begin{df}\label{def:de}
If $q \in S(\mathfrak{h})$ we denote by 
\index{dqk@$D_q^\pm(k)$}$D_q^\pm(k)$ the differential 
operator that coincides with $q(T_\xi(k))$ on $\mathbf{C}[H]^{\pm W}$.
\end{df}

\begin{lm}\label{lm:equiv}
We put
\index{pikpm@$\pi^\pm(k)$}
\[
\pi^\pm(k)
 = \prod_{\alpha \in R_+} (\alpha^\vee \pm k_\alpha) 
\in S(\mathfrak{h}) \subset \mathbf{H},\] 
and denote by \index{idemp@$\varepsilon^\pm$}$\varepsilon^\pm$
the idempotents in $\mathbf{C}[W]$ corresponding 
to the trivial representation ($\varepsilon^+$) and the sign representation 
($\varepsilon^-$) respectively. Then
\begin{itemize} 
\item[(a)]
$
\varepsilon^\mp \cdot \pi^\pm(k)\cdot \varepsilon^\pm = \pi^\pm(k)\cdot \varepsilon^\pm.
$   
\item[(b)]
$
\varepsilon^\pm \cdot \mathbf{H}(k)\cdot \varepsilon^\pm = Z(\mathbf{H}(k))\cdot \varepsilon^\pm. 
$
The map $Z(\mathbf{H}(k))\to Z(\mathbf{H}(k))\cdot \varepsilon^\pm$, $z\to 
z\cdot \varepsilon^\pm$ is an isomorphism of commutative algebras, and the map 
${Rad}^\pm:\mathbf{H}(k) \to Z(\mathbf{H}(k))$ defined by $\varepsilon^\pm \cdot h\cdot 
\varepsilon^\pm={Rad}^\pm(h)
\cdot \varepsilon^\pm$ respects the filtering by 
degree.

\item[(c)]
$
\varepsilon^\mp \cdot \mathbf{H}(k)\cdot \varepsilon^\pm = Z(\mathbf{H}(k))\pi^\pm(k)\cdot
\varepsilon^\pm
$
The map $Z(\mathbf{H}(k))\to Z(\mathbf{H}(k))\pi^\pm(k)\cdot \varepsilon^\pm$, $z\to 
z\pi^\pm(k)\cdot \varepsilon^\pm$ is a linear isomorphism, and the map 
${}^\pm{Rad}:\mathbf{H}(k) \to Z(\mathbf{H}(k))\pi^\pm(k)$ 
defined by $\varepsilon^\mp \cdot h\cdot 
\varepsilon^\pm={}^\pm{Rad}(h)
\cdot \varepsilon^\pm$ respects the filtering by 
degree.
\end{itemize}
\end{lm}
\pf 
To prove (a) it is enough to show that for all simple reflections $r_i$, 
\[
(r_i\cdot\pi^\pm(k)+\pi^\pm(k)\cdot r_i)\cdot\varepsilon^\pm=0.
\]
This follows from Definition~\ref{def:hecke} (4).
As to (b), first observe that it is enough to show that for all $p\in S(\mathfrak{h})$, 
$\varepsilon^\pm \cdot p\cdot \varepsilon^\pm \in Z(\mathbf{H}(k))\cdot \varepsilon^\pm$.
Using formula \ref{eqn:wxi} and Lemma \ref{lm:cent} this is clear, 
by induction on the degree of $p$. The remaining statements follow trivially from 
this proof. Essentially 
the same arguments, combined with (a), proves (c). 
\qed 
\begin{df}\index{fundamental shift operators}
The {\it fundamental shift operators} 
\index{gk@$G_\pm(k)$}$G_\pm(k)$ are defined by 
\[
G_+(k) = \Delta^{-1} D_{\pi^+(k)}^+(k),
\]
and 
\[
G_-(k+1) = D_{\pi^-(k)}^-(k)\Delta.
\]
\end{df}
The shift principle is equivalent with the following action of the shift  
operators on Jacobi polynomials:
\begin{thm}\label{thm:shiftrels}
 We have the following shift relations ($\lambda\in P_+$):
\[
G_+(k)P(\lambda,k)=\prod_{\alpha\in R_+}(k_\alpha-(\lambda+\rho(k)(\alpha^\vee))P(\lambda-\delta,
k+1)
\]
and 
\[
G_-(k+1)P(\lambda,k+1)=\prod_{\alpha\in R_+}(k_\alpha+(\lambda+\delta+\rho(k)(\alpha^\vee))P(\lambda+\delta,
k)
\]
\end{thm}
\pf
Both relations are proved in the same manner. Let us do the first one.
By Lemma \ref{lm:equiv} it is clear that 
\[
D_{\pi^+(k)}^+(k)P(\lambda,k)=c\cdot P^-(\lambda,k)
\]
for some constant $c$. To compute this constant one has to recall that the 
Dunkl operators are triangular with respect to the ordering $\triangleleft$.
With respect to this ordering, the highest order term in the expansion of 
$P(\lambda,k)$ is $e^{w_0\lambda}$, and the highest order term of 
$P^-(\lambda,k)$ is $\varepsilon(w_0)e^{w_0\lambda}$. Using Corollary 
\ref{cor:ev} and the shift principle it is now straightforward to verify 
the asserted relation. 
\qed
We collect some basic properties of the shift operators in the 
following theorem.
\begin{thm}\label{thm:shiftprop}
$(a)$ $G_\pm(k)$ transforms $\mathbf{C}[\h]^W$ to $\mathbf{C}[\h]^W$  \\
$(b)$ For all 
$f,g \in \mathbf{C}[H]$, $(G_+(k)f,g)_{k+1} = (f,G_-(k+1)g)_k $    \\
$(c)$ For all $p \in S(\mathfrak{h})^W$, 
$D_p(k\pm 1) G_\pm(k) = G_\pm (k)D_p(k)$ \\
$(d)$ For any $W$-invariant holomorphic germ $f$ at $x=e$, we have 
\[
(G_-(k+1)f)(e) = 
\frac{\tilde{c}(\rho(k),k)}{\tilde{c}(\rho(k+1),k+1)} f(e).
\]
\end{thm}
\pf
(a) In the case of $G_+(k)$ this is immediate from Remark 2.5, and 
in the case of $G_-(k)$ we use \ref{lm:equiv} and 
the divisibility of $W$-skew Laurent polynomials 
by $\Delta$. 

(b)
From the definitions and the symmetry of the Dunkl-Cherednik operators 
with respect to the inner product $(\cdot,\cdot)_k$, we see that one has to 
verify (in the terminology of Lemma \ref{lm:equiv} (c)) that ${}^-Rad(\pi^+(k))
={}^-Rad(\pi^-(k))$. This is true because Lemma \ref{lm:equiv} (c) implies that 
 ${}^-Rad$ kills polynomials with degree 
lower than $|R_+|$.

(c) This is an immediate consequence of Theorem \ref{thm:shiftrels}.

(d) By powerseries expansion at $e$ it is easy to see that 
\begin{equation}\label{eqn:ate}
(G_-(k+1)f)(e)=c\cdot f(e)
\end{equation}
for a some constant $c$. When we apply this to the function 
$f=1=P(0,k+1)$ and use Theorem \ref{thm:shiftrels} we find that 
\begin{equation}\label{eqn:valc}
c=\prod_{\alpha\in R_+}(k_\alpha+(\delta+\rho(k))(\alpha^\vee))P(\delta,k,e)
\end{equation}
Taking $f=P(\lambda,k+1)$ in \ref{eqn:ate} we now obtain 
\begin{align*}
P(\lambda,k+1,e)&P(\delta,k,e)\prod_{\alpha\in R_+}(k_\alpha+(\delta+\rho(k))(\alpha^\vee))=\\
&P(\lambda+\delta,k,e)\prod_{\alpha\in R_+}(k_\alpha+(\lambda+\delta+\rho(k))(\alpha^\vee)).
\end{align*}
This is a recursive formula for $P(\lambda,k,e)$, that can be 
solved starting from $P(\lambda,0,e)=|W_\lambda|$.
This quickly leads to the formula 
\begin{equation}\label{eqn:valpp}
P(\lambda,k,e) = 
\frac{ \tilde{c}(\rho(k),k)}{\tilde{c}(\lambda+\rho(k),k)}.
\end{equation}
Now the constant $c$ from equation \ref{eqn:ate} follows from 
\ref{eqn:valc} and \ref{eqn:valpp}.  
\qed 
\begin{cor}{\rm (of proof)}\label{cor:comp}
The value of $P(\lambda,k,e)$ (see equation \ref{eqn:valpp}) can be computed by the 
use of the shift operators. Likewise we can compute the square norms of the Jacobi 
polynomial $P(\lambda,k)$ with respect to $(\cdot,\cdot)_k$ by a recursion relation 
that follows from Theorem \ref{thm:shiftrels} and Theorem \ref{thm:shiftprop} (b). Details 
are left to the reader (see \cite{Op02}).
\end{cor}

\begin{rem}
Obviously the square norms and special values of 
the Jacobi polynomials obtained in Corollary \ref{cor:comp} 
could have been obtained immediately from   
Theorem~\ref{thm:norm}. However, the converse is also true, 
up to some algebraic manipulations in $\mathbf{H}(k)$ 
(see \cite{Op}, Section 5). In other words, with respect to 
the results of Theorem~\ref{thm:norm} both the method of affine 
intertwiners (Section 4) and the method of shift operators 
are simple and effective. (This remark is true in the cases of 
the Macdonald and Koornwinder orthogonal polynomials as well.) 
However, because we use division by $\Delta$ in the generalized 
Weyl character formula, the shift operators are not suitable for 
proving combinatorial formulae, or the positivity and integrality 
conjectures.   
\end{rem}

\subsection{Another reflection symmetry and application}

The operators $D_p(k)$ have another symmetry in the parameter $k$ that 
gives a direct relation between the two shift operators $G_-$ and $G_+$.
This has an important application because it gives a proof of the 
conjecture by Yano and Sekiguchi concerning the explicit form of the
$b$-function  for the discriminant of a crystallographic reflection group.
\begin{theorem}{\rm (see Proposition 2.2 of \cite{HO1})} Let $1-k\in K$ be
defined by $(1-k)_\alpha=1-k_\alpha$. Then we have:
\[
D_p(1-k)=\delta_{k-1/2}\circ D_p(k) \circ \delta_{1/2-k}.
\]
\end{theorem}
\pf (Sketch)
When $p_2=\sum x_i^2$ this is a direct computation using the
explicit  formula in Example \ref{ex:L} for
$D_{p_2}(k)=L(k)+(\rho(k),\rho(k))$. It is not difficult and 
standard to see that an operator $D$ that commutes with 
$D_{p_2}(1-k)$, and that has an asymptotic expansion as in 
\ref{eqn:asexp}, is determined by its image $p=\gamma(D)$ under the
Harish-Chandra homomorphism (see also \ref{eqn:HCH}). Therefore the 
conjugation formula holds for all $p\in S(\mathfrak{h})^W$. \qed

By a similar argument one proves the following consequence:
\begin{cor}
\[
G_+(-1/2-k)\circ\delta_{k+1}=\delta_k\circ G_-(3/2+k)
\]
\end{cor} 
Now apply this identity to the constant function $1$, and take 
the lowest homogeneous part of the identity thus obtained. Use 
\ref{thm:shiftprop}(d). This gives: 
\begin{cor} Take $k_\alpha=k\ \forall\alpha\in R$. 
Let $D$ be the lowest
homogeneous part of  $G_+(-1/2-k)$ at the unit element of $H$. Let 
\[
\pi=\prod_{\alpha\in R_+}\alpha^2
\]
be the discriminant of the reflection group $W$.
Then
\[
D\pi^{k+1}=|W|\prod_{i=1}^n\prod_{j=1}^{d_i-1}(d_i(k+1/2)+j)\pi^k.
\]
where $d_1,\dots,d_n$ are the primitive degrees of $W$.
\end{cor}
From this formula it is easy to compute the $b$-function 
of $\pi$. The result was conjectured by Yano and Sekiguchi in 
\cite{YS}.
\begin{theorem}{\rm (\cite{Op02}, Theorem 7.1)}
The $b$-function of the discriminant $\pi$ is given by:
\[
b(s)=\prod_{i=1}^n\prod_{j=1}^{d_i-1}(s+1/2+\frac{j}{d_i}).
\]
\end{theorem}
\begin{rem} We have introduced two shift operators $G_\pm$ in 
this section, associated to the sign character of $W$. 
In fact one can associate a raising and a 
lowering operator to each linear character of $W$. For the 
purpose of this section we did not need this construction so 
we have skipped it. The interested reader is advised to consult 
\cite{Op02} and \cite{HS} for the properties of these shift 
operators. 
\end{rem}
\newpage
%
\section{Away from polynomials }

This section is a review of 
the hypergeometric function for root systems, 
which is a $k$-deformation of 
the elementary spherical function on symmetric spaces.
This function was introduced and studied by Heckman and 
Opdam in  \cite{HO1} and a series of subsequent papers.
An introduction to 
the hypergeometric system and the hypergeometric 
function is \cite[Part I]{HS}, where one can find further 
references.

In the previous section,
we have introduced the differential operator 
\index{dpk@$D_p(k)$} $D_p(k) = D_p^+(k)$ for $p \in S({\mathfrak h})^W$,
which maps ${\bf C}[H]^W$  to itself.
By Chevalley's theorem
${\bf C}[H]^W \cong {\bf C}[z_1,z_2,\dots,z_n]$
with \index{zi@$z_i$} $z_i = \sum_{\mu \in W \mu_i} e^{\mu}$,
so we have a system of commuting partial differential operators
on the affine space $W \backslash H$.
We want to study the general eigenvalue problem for these
operators.
We have seen that
when we want polynomial eigenfunctions $\varphi\in {\bf C}[H]^W$, 
we are forced to take the eigenvalue $\lambda \in \mathfrak h^*$ in the 
system
\[
D_p(k) \varphi = p(\lambda) \varphi, \qquad \forall p \in S(\mathfrak h)^W
\]
equal to $\mu + \rho(k)$ for some $\mu \in P_+$.
This means that the eigenvalue has to satisfy 
a certain integrality condition in this stiuation.
However, for values of $\lambda$ that are not integral in this sense, 
we can still find germs of holomorphic solutions at 
any point $h\in H$. The most elementary case is the 
case where $h$ is regular for the action of $W$. 
We will see in the next subsection that in this case 
the space of germs of holomorphic solutions has dimension 
$|W|$. For generic parameters we can give a basis of 
series solutions which are convergent in an open neighbourhood 
of $A_+$, and which behave asymptotically free (the Harish-Chandra 
series). 

The important conclusion at this point is that the sheaf of 
germs of holomorphic solutions of these equations 
(equations \ref{eqn:hgs}) is a local system of rank $|W|$ 
on the regular $W$ orbit space of $H$. A further understanding 
of the equations \ref{eqn:hgs} is obtained from the investigation 
of the monodromy of the local system, in subsection 6.2.

\subsection{Harish-Chandra series} 
We denote the set of regular elements by
\index{hreg@$H^{\text{reg}}$}
\[
H^{\text{reg}}
=\{ h \in H \,;\, \Delta^2(h)\neq 0 \}.
\]
We choose a base point $z \in  W \backslash H^{\text{reg}}$
with a representative $h \in H^{\text{reg}}$.
By definition,
the germ ${\cal O}_z$ of holomorphic functions at $z$
is the germ $\cong {\cal O}_{W h}^W$ 
of $W$-invariant holomorphic functions on $W h$.
Remark that ${\cal O}_{W h} = \oplus_{w\in W} {\cal O}_{w h}$.
\begin{dfn}\index{hypergeometric system}
The {\em hypergeometric system} of differential equations
at $z \in W \backslash H^{\text{reg}}$
with a spectral parameter $\lambda \in \mathfrak{h}^*$
is the system of differential equations
\begin{equation}\label{eqn:hgs}
D_p(k) \varphi = p(\lambda) \varphi, 
\qquad p \in S(\mathfrak h)^W
\end{equation}
for an unknown function $\varphi \in {\cal O}_z \cong {\cal O}_{W h}^W$.
\end{dfn}
We denote the set of solutions for this system by 
\index{slambdakw@$S(\lambda,k)^W$}
\[
S(\lambda,k)^W
= \{ \varphi \in {\cal O}_{W h}^W 
  \,;\, D_p(k) \varphi = p(\lambda) \varphi, \,
 p \in S(\mathfrak h)^W \}.
\]
\begin{example}\label{ex:L}
Let $\xi_1,\dots,\xi_n$ be an orthonormal basis of $\mathfrak{a}$.
Then  $p = \sum \xi_i^2$ is a $W$-invariant quadratic,
and the corresponding differential operator is
\[
D_p(k) = L(k) + (\rho(k),\rho(k)),
\]
where
\index{lk@$L(k)$} 
\[
L(k)
= \sum_{i=1}^n \partial_{\xi_i}^2 
+ \sum_{\alpha \in R_+} \frac12 k_\alpha 
\frac{1+e^{-\alpha}}{1-e^{-\alpha}} (\alpha,\alpha) 
\partial_{\alpha^\vee}.
\]
Let $\mathfrak{g}$ be a real semisimple Lie algebra
with Cartan decomposition ${\mathfrak g} = {\mathfrak k} \oplus {\mathfrak p}$
and $\mathfrak a \subset \mathfrak p$ a maximal abelian subspace,
and $\Sigma = \Sigma(\mathfrak g, \mathfrak a)$  the
restricted root system with root labels 
$m_\alpha = \dim(\mathfrak g^{\alpha})$. 
Then the radial part of the Laplace-Beltrami operator on $G/K$
with respect to left action of $K$ equals $L(k)$, 
if we identify 
$R$ with $2\Sigma$ and
$k_{2\alpha}=\frac12 m_\alpha$.
So (\ref{eqn:hgs}) 
becomes the system of differential equations for the 
\index{elementary spherical function}
elementary spherical function $\varphi_\lambda$ 
restricted to $A$.
\end{example}
\begin{example}\label{ex:rk1}
Let us consider the rank $1$ case,
and in order to be even more convincing,
we do the non-reduced case BC$_1$, 
$R = \{ \pm \alpha, \pm 2 \alpha \}$.
Let us introduce notation.
$H = {\bf C}^\times$,
${\bf C}[H] = {\bf C}[y,y^{-1}]$, with
$y=e^\alpha$ ;
If $\xi = (2\alpha)^\vee$, then $Q^\vee = P^\vee$
is generated by $\xi$,
and $\partial_\xi = \theta = y \frac{d}{dy}$.
Normalize $|\xi|=1$.
We set $\lambda=\lambda(\xi)$, $k_1=k_\alpha$, $k_2=k_{2\alpha}$.
Now (\ref{eqn:hgs}) becomes
\[
\left\{ \theta^2 
+ \left( k_1 \frac{1+y^{-1}}{1-y^{-1}} 
     + 2 k_2 \frac{1+y^{-2}}{1-y^{-2}} \right) \theta 
+ \left( (\frac12 k_1 + k_2)^2 - \lambda^2 \right) \right\} \varphi = 0.
\]
Let $z = \frac12 - \frac14(y+y^{-1})$ be 
a coordinate on $W \backslash H$,
then this becomes
\[
\left\{ z(1-z) \frac{d^2}{dz^2} 
+ (c - (1+a+b) z)\frac{d}{dz}
-ab \right\} \varphi = 0
\]
with $a=\lambda + \frac12 k_1 + k_2$,
$b=-\lambda + \frac12 k_1 + k_2$, $c = \frac12 + k_1 + k_2$.
\end{example}
%
To understand system (\ref{eqn:hgs}),
we first consider the easiest examples of solutions,
the asymptotically free solutions on $A_+$ 
(also called the Harish-Chandra series). 
\index{Harish-Chandra series}

The crucial point is the observation that the equations 
themselves have an asymptotic expansion as follows.
\begin{lemma}\label{lemma:asexp}
For any $p\in S({\mathfrak h})^W={\bf C}[{\mathfrak h^*}]^W$ one has an  
asymptotic expansion of the following kind on $A_+$:
\begin{equation}\label{eqn:asexp}
D_p(k)=\partial(p(\cdot +\rho(k)))+\sum_{\kappa\in
Q_-\backslash\{0\}} e^\kappa\partial(p_\kappa)
\end{equation}
where $p_\kappa\in {\bf C}[{\mathfrak h^*}]$ has lower degree than $p$.
More generally, for any $p$ in $S({\mathfrak h})={\bf C}[{\mathfrak h^*}]$ 
and $w_0\in W$ the longest element of $W$, we have the following asymptotic
expansion on $A_+$ (compare with 
\cite{Op}, Lemma 6.4):
\[
w_0D_{p^{w_0}}(k)w_0=\partial(p(\cdot + \rho(k)))+\sum_{\kappa\in
Q_-\backslash\{0\}} e^\kappa\partial(p_\kappa)
\]
\end{lemma}
\pf
We prove the second asymptotic formula, 
by induction on the degree of $p$. 
Let $p$ be of the form 
$p=\xi q$ with $\xi\in\mathfrak h$ 
and let $w\in W$. 
Then 
\begin{gather}
\begin{split}\label{eq:asymp}
w^{-1}D_{(\xi q)^{w}}&(k)
w=(\partial_\xi-w^{-1}\rho(k)(\xi))w^{-1}D_{q^{w}}(k)w+ \\
&\sum_{\alpha\in
R_+}k_\alpha\alpha(\xi)\frac{1}{1-e^{-w^{-1}\alpha}}
(w^{-1}D_{q^{w}}(k)w-w^{-1}r_\alpha D_{q^{w}}(k)
r_\alpha w)
\end{split}
\end{gather}
(just check that the right hand side is a differential 
operator that restricts to $w^{-1}(p^{w}(T))w$ 
on $W$-invariant functions).
From \ref{eq:asymp} it follows by induction that $w^{-1}D_{p^w}w$ 
has an asymptotic expansion on $A_+$ of the form:
\begin{equation}\label{eq:genas}
\sum_{\kappa\in
Q_-}e^\kappa\partial(p_\kappa)
\end{equation}
with $\text{deg}(p_\kappa)\leq\text{deg}(p)$, 
with equality if and only if $\kappa=0$.
In the special 
case where $w=w_0$ we want to prove that 
$p_0(\lambda)=p(\lambda+\rho(k))$. 
Observe that in this special case none of the 
terms of the second line of \ref{eq:asymp} contribute to 
the leading term (using \ref{eq:genas}).
Hence the result follows from \ref{eq:asymp} by induction on 
the degree. 
\qed

Substitute a formal series 
\[
\varphi=\sum_{\nu \le \mu} c_\nu e^{\nu},
\qquad c_\mu=1
\]
into (\ref{eqn:hgs}). By lemma \ref{lemma:asexp} we 
obtain the following indicial equation for the leading 
exponent: 
\begin{equation}\label{eq:ind}
p(\mu+\rho(k)) = p(\lambda),
\qquad p \in S(\mathfrak{h})^W.
\end{equation}
This means that 
\[
\lambda \in W (\mu + \rho(k)).
\]
We put $\lambda = \mu + \rho(k)$, and put $c_\nu =
\Gamma_\kappa(\lambda,k)$ 
\index{gammakappa@$\Gamma_\kappa(\lambda,k)$}
if $\kappa = \nu - \mu \in Q_-$.
Just using the explicit second order operator $L(k)$ 
we arrive at the following recurrence relations.
\begin{equation}\label{(**)}
-(2\lambda+\kappa,\kappa) \Gamma_\kappa(\lambda,k)
=
2 \sum_{\alpha>0} k_\alpha \sum_{j\ge1}
(\lambda-\rho(k) + \kappa + j \alpha,\alpha)
\Gamma_{\kappa+j \alpha}(\lambda,k)
\end{equation}
These have a unique solution if we fix $\Gamma_0(\lambda,k)=1$, 
and then the coefficients $\Gamma_\kappa(\lambda,k)$
are rational, with poles possibly at the hyperplanes
$H_{\kappa'}$ for some  $\kappa'<0$, where 
\index{hkappa@$H_\kappa$}
\begin{equation}\label{eq:hk}
H_\kappa =
\{\lambda \in {\mathfrak h} \,;\,
(2\lambda + \kappa,\kappa) = 0 \}.
\end{equation}
Next we want to show that the eigenfunctions of the 
second order equation which we have just constructed,  
are in fact solutions of all 
the  equations \ref{eqn:hgs}. The following well known 
and beautiful argument is due to Harish-Chandra. 
The uniqueness of the asymptotic solution, 
combined with the lemma \ref{lemma:asexp}
and the commutativity of the operators 
$\{D_p\,;\,p\in S(\mathfrak{h})^W\}$
imply that  
\index{phil@$\Phi(\lambda,k)$}
\[
\Phi(\lambda,k)
 = \sum_{\kappa \in Q_-} \Gamma_\kappa(\lambda,k)
e^{\lambda - \rho(k) + \kappa}, \quad
\Gamma_0(\lambda,k)=1
\]
is a joint eigenfunction of the commuting family of 
differential operators $\{D_p\,;\,p\in S(\mathfrak{h})^W\}$. 
It is easy to find the eigenvalues by considering the 
leading exponents, taking lemma \ref{lemma:asexp} into 
consideration. We find that 
\[
D_p(k)\Phi(\lambda,k)=p(\lambda)\Phi(\lambda,k)
\]
In other words, we have indeed constructed formal series 
solutions of \ref{eqn:hgs}.
In this context one traditionally writes 
\begin{equation}\label{eqn:HCH}
p(\lambda)=\gamma(D_p(k))(\lambda),
\end{equation} 
and then one 
calls $\gamma$ the ``Harish-Chandra homomorphism''.

The series $\Phi(\lambda,k)$ converges on 
\index{aplus@$A_+$}
\[A_+
 = \{ a \in A \,;\, 
a^\alpha = e^{\alpha}(a) >1,\quad \forall \alpha >0 \}.
\]
as one easily verifies using the defining recurrence relations.
 
As we have seen in the descriptions above, 
there are possibly singularities in the parameter space $\mathfrak{h}^*\times
K$ of our series solutions $\Phi(\lambda,k)$. These are simple poles along the
hyperplanes $H_\kappa$ as defined in \ref{eq:hk}. However, 
the actual set of
poles $\Phi(\lambda,k)$ turns out to be a much smaller subset of
hyperplanes:
\begin{lemma}\label{lem:pole}
The (apparent) simple pole of $\Phi(\lambda,k)$ (as a function of 
$\lambda$!) along
$H_\kappa$  is removable unless $\kappa=n\alpha$ for some $n\in\Z_-$ and
$\alpha\in R_+$. If $\kappa=n\alpha$ then the residue of $\Phi(\lambda,k)$ at 
$H_\kappa$ is a constant multiple of $\Phi(r_\alpha(\lambda),k)$.
\end{lemma}
\pf
From the recurrence relations it is easy to see that the residue 
of $\Phi(\lambda,k)$ at $H_\kappa$ is a constant multiple of 
$\Phi(\lambda+\kappa,k)$. Suppose it is nonzero. Then by the indicial 
equation \ref{eq:ind}, the leading exponent $\lambda+\kappa$ of the residue 
must be of the form $w\lambda$ for some $w\in W$, and this must hold 
for all $\lambda\in H_\kappa$. Hence $w=r_\alpha$ for some $\alpha\in R$, 
and $\kappa=n\alpha$ for some $n\in \Z$. It is obvious that $\kappa$ has to be 
negative in the dominance ordering.
\qed

The equation that defines $H_{n\alpha}$ can be rewritten as 
\[
\lambda(\alpha^\vee)+n=0.
\]
We now change the notation for this hyperplane  
to $H_{n,\alpha}$, so as to also include the case $n=0$ of 
the hyperplane perpendicular to the root $\alpha$.
We will call $\lambda$ generic if 
\begin{equation}\label{eq:generic}
\lambda\not\in \cup_{n\in\Z,\alpha\in R} H_{n,\alpha}.
\end{equation}
\begin{remark} Notice that the set of generic parameters 
is precisely the set of regular points for the action of 
the affine Weyl group introduced in Section 3. There is a natural
action of the affine  Weyl group on the space of nonsymmetric eigenfunctions
of the Dunkl operators $T_\xi$, via the intertwiners of
Section 4. The relation between such nonsymmetric eigenfunctions and our  
space of solutions of \ref{eqn:hgs} is the subject of the next section.
\end{remark}
If $\lambda$ is generic
then, by Lemma \ref{lem:pole}, the dimension  of the solution space for the
eigenfunction equations \ref{eqn:hgs} on $A_+$ is  at least equal to $|W|$.
The next theorem tells us that this  is in fact an equality which holds for
any $\lambda$, and moreover  that this is the dimension of the solution space
of these equations  in the space of holomorphic germs at any regular point of
$H$.

\begin{theorem}\label{thm:rank}
System {\rm (\ref{eqn:hgs})} is holonomic of rank $|W|$.
If $\lambda\in\h^*$ is generic 
then 
$\{\Phi(w\lambda,k;\cdot)\,;\,w\in W\}$ 
forms a basis of the solution space.
\end{theorem}
\pf 
For any homogeneous $p \in S(\mathfrak h)^W$,
\[
D_p(k) = \partial(p) + \mbox{(lower order terms)}.
\]
Then in the left ideal generated by $D_p(k) - p(\lambda)$,
we have operators of the form
\[
\partial(q) + \mbox{(lower order terms)},
\qquad \forall q \in S(\mathfrak h)S(\mathfrak h)^W_+,
\]
where $S(\mathfrak h)^W_+$ denotes the space of the elements of
$S(\mathfrak{h})^W$ without constant term. 
Hence the left ${\cal O}_z$-module
\[
{\cal D}_z /  
\sum _{p \in S(\mathfrak h)^W} {\cal D}_z ( D_p(k) - p(\lambda) )
\]
is generated by the operators 
\[
\partial(q),
\qquad \mbox{with $q \in S({\mathfrak h})$, $W$-harmonic polynomials}.
\]
Then the holonomic rank at the base point $z$ 
is less than or equal to $|W|$.
Conversely, we found, generically,
the linearly independent asymptotically free solutions 
$\Phi(w \lambda, k; \cdot)$.
Combining these, we conclude that
the holonomic rank equals $|W|$ generically.

A more precise version of this argument shows that
$(\partial(q))$ ($q\in S({\mathfrak h})$: harmonic) always gives 
an ${\cal O}_z$-basis for the ${\cal D}_z$-module,
independent of the parameter choice (see \cite{HO1} or 
\cite{HS}). This point will also become quite clear 
in section 7, when we study the relation with between 
\ref{eqn:hgs} and the KZ connection.\qed

\subsection{Monodromy} 
We need to understand the {\it monodromy action}
of $\pi_1(W \backslash H^{\text{reg}},z_0)$
on the solution space of (\ref{eqn:hgs}).
Take a base point $x_0 \in A_+ \subset A^{\text{reg}}$
such that $z_0 = \overline{x_0}$.
For each simple reflection $r_i$ 
we consider an element $l_i$ in $\pi_1(W \backslash H^{\text{reg}},z_0)$
defined as follows: $l_i$ can be represented by a path from $x_0$ to
$r_i(x_0)$ which we can take arbitrarily close to the "straight" line 
segment between these two end points, but near the wall $a^{\alpha_i}=1$ 
we replace a subsegment that intersects the wall by a half circle going around
the wall in positive direction. 

For each 
$v \in Q^\vee$ we define the closed loop
\index{lv@$l_v$} $l_v$   by 
\[
l_v(t) = x_0 \exp(2\pi \sqrt{-1} t v)\quad(t \in [0,1]).
\]
Given $\varphi$, a local solution at $x_0$ of (\ref{eqn:hgs}),
we denote $T_i \varphi$ 
for the solution obtained by continuing $\varphi$ 
analytically along the path $l_i$, and composing the result with
$r_i$, and we denote $T_v \varphi$ for the continuation of $\varphi$
along the loop $l_v$.

\medskip
System (\ref{eqn:hgs}) has regular singularities at infinity 
and and also along the walls. Moreover the structure of 
the fundamental group 
$\pi_1(W \backslash H^{\text{reg}},z_0)$ allows the method of 
rank one reduction, which enables us to compute the connection 
formula  for $\{\Phi(w\lambda,k;\cdot)\,;\,w\in W\}$ explicitly 
in terms of the $c$-function:
\begin{theorem}{\rm (Looijenga, v.d.Lek, Heckman-Opdam)}\label{thm:mono}
Assume that $\lambda\in\h^*$ satisfies condition \ref{eq:generic}. 
\begin{itemize}
\item[{\rm (a)}] 
Put $T_0 = T_{\theta^\vee} T_{i_1} \cdots T_{i_k}$
with $r_{i_1}\cdots r_{i_k}$ a reduced expression for $r_{\theta^\vee}$.
This is independent of the reduced expression,
and $T_0$, $T_1$, \dots, $T_n$ satisfy the 
braid relations of $W^a$.
These operators generate all monodromy on 
$W \backslash H^{\text{\rm reg}}$ (in other words, the 
corresponding elements of $\pi_1(W \backslash H^{\text{reg}},z_0)$
form a set of generators).
\item[{\rm (b)}] 
$(T_i-1)(T_i+q_i)=0$ for all $i=0,1,\dots,n$,
with $q_i = e^{-2\pi\sqrt{-1} k_i}$.
\item[{\rm (c)}] 
$T_v \Phi(\lambda,k) = e^{2 \pi \sqrt{-1} (\lambda - \rho(k))(v)}
\Phi(\lambda,k)$.
\item[{\rm (d)}]
$\tilde{c}(\lambda,k) \Phi(\lambda,k)
+ \tilde{c}(r_i\lambda,k) \Phi(r_i\lambda,k)$
is fixed for $T_i$ ($i=1,\dots,n$).
\item[{\rm (e)}]
$\tilde{c}(-r_i\lambda,1-k) \Phi(\lambda,k)
+ \tilde{c}(-\lambda,1-k) \Phi(r_i\lambda,k)$
has eigenvalue $-q_i$ with respect to $T_i$ 
($i=1,\dots,n$).
\end{itemize}
\end{theorem}
\pf 
As indicated, these results come from various sources; we refer to 
\cite[Part I, Lecture 4]{HS} 
for more details and references. 
\begin{itemize}
\item[(a)]
is  known from the work of Looijenga and v.d.Lek 
on the fundamental group 
$\pi_1(W \backslash H^{\text{reg}}, x_0)$,
and is a nontrivial result.
\item[(b)] 
follows from (d) and (e).
\item[(c)]
is trivial.
\item[(d)] and (e) form the heart of the matter. 
The proof is not difficult, and reduces to the rank one case. 
Let us sketch the idea of the proof. 
From the braid relations (a) 
it follows that if $v \in Q^\vee$ such that
$\alpha_i(v) = 0$,
then \index{ti@$T_i$} $T_i$ and \index{tv@$T_v$} $T_v$ commute
(already in the fundamental group).
Hence by (c) we see that, for generic $\lambda$, 
span$(\Phi(\lambda,k), \Phi(r_i\lambda,k))$
is closed for $T_i$.
Now one takes limiting values of
\[
e^{-\lambda+\rho(k)} \Phi(\lambda,k,b\cdot\exp(t\alpha^\vee_i))
\]
when $b \to \infty$ in the wall $b^{\alpha_i}=1$.
%
%
The resulting limits are formal series solutions
(asymptotically free at $\infty$) 
of Example \ref{ex:rk1},
and here the monodromy of such series is explicitly known. 
For the precise argument, see \cite[Theorem~6.7]{HO1}, 
\cite[Theorem~1.1]{Hec1}, and \cite[Part I, Lecture 1,
Section~4.3]{HS}. 
\end{itemize}
\qed 
Motivated by these facts, we define
the affine Hecke algebra \index{har@${\bf H}^{\mathrm{aff}}(R_+, q_i)$} 
${\bf H}^{\mathrm{aff}}(R_+, q_i)$
generated by $T_i$'s and $T_v$'s
with the relations (a) and (b) in Theorem \ref{thm:mono}
This algebra contains two important subalgebras;
the finite dimensional Hecke algebra ${\bf H}(R_+,q_i)=\langle T_i
\rangle_{i=1}^n$ (describing the monodromy locally at the 
unit element of $H$), and the group algebra
${\bf C}[Q^\vee]=\langle \theta_v \rangle_{v \in
Q^\vee}$, where $\theta_v$ is defined by
$\theta_v=e^{2\pi\sqrt{-1}\rho(k)(v)}T_v$ (describing the monodromy ``at
infinity'' in $A_+$).  As a vector space, the algebra ${\bf H}^{\mathrm{aff}}(R_+, q_i)$
is naturally isomorphic to the tensor product of these two algebras: \[ {\bf
H}^{\mathrm{aff}}(R_+,q_i) \cong {\bf H}(R_+,q_i) \otimes {\bf C}[Q^\vee].\]
The relations between the $T_i$ and the $\theta_v$ are given by Lusztig's 
formula:
\begin{equation}\label{eqn:lus}
T_i\theta_v-\theta_{r_iv}T_i=(q_i-1)
\left(\frac{\theta_v-\theta_{r_iv}}{1-\theta_{-\alpha_i^\vee}}\right)
\end{equation}
\begin{cor}\label{cor:monoind}
The monodromy is, for generic parameters, 
equal to the representation 
\[
\mbox{\rm Ind}_{{\bf C}[Q^\vee]}^{{\bf H}^{\mathrm{aff}}(R_+,q_i)}
e^{2\pi\sqrt{-1}(\lambda)}
\]
Here we consider $e^{2\pi\sqrt{-1}(\lambda)}$
as a character of ${\bf C}[Q^\vee]$.
\end{cor}

\begin{remark}\label{rem:regsing}
At this point it is natural to invoke the result 
that the holonomic system of differential equations \ref{eqn:hgs} has regular 
singularities, both at the ``hyperplanes'' $e^\alpha=1$ in $H$ and 
``at infinity'' when we consider the torus $H$ as a quasi-projective 
variety (for instance via an embedding in a projective toric variety).
These facts have simple proofs which will be given 
in section 7, when we study the equivalence of 
\ref{eqn:hgs} and the KZ connection. The point is that the 
KZ connection visibly meets these regularity requirements.
\end{remark} 
\begin{cor}\label{cor:cont} Let $\lambda$ be generic.
The linear combination $\tilde{c}(\lambda,k) \Phi(\lambda,k)
+ \tilde{c}(r_i\lambda,k) \Phi(r_i\lambda,k)$ 
as mentioned in  Theorem \ref{thm:mono}(d)
extends holomorphically in a neighbourhood of 
$int(\overline{A_+\cup r_i(A_+)})$, and is $r_i$ invariant. 
Hence for generic $\lambda$,
the function (for $\tilde{c}$: see\ref{eqn:ctilde}) 
\index{ft@$\tilde{F}(\lambda,k;a)$}
\[
\tilde{F}(\lambda,k;a)
= \sum_{w \in W} \tilde{c}(w\lambda,k) \Phi(w\lambda,k;a)
\]
extends holomorphically from $A_+$ to a tubular
neighbourhood of $A$ in $H$,
and is $W$-invariant there.
\end{cor}
\pf
The linear combination of Harish-Chandra series 
under consideration has no monodromy with respect to 
$l_i$ by \ref{thm:mono}, which means that it extends to a 
$r_i$-invariant holomorphic function on an open set of the
form $U\cdot int(\overline{A_+\cup r_i(A_+)})\backslash  \{e^{\alpha_i}=1\}$
where $e\in U$, $U\subset T$ open and connected.
By Remark  \ref{rem:regsing} this function has moderate growth towards 
$\{e^{\alpha_i}=1\}$, hence it will extend meromorphically 
to $int(\overline{A_+\cup r_i(A_+)})$. Let us denote its pole order along 
$\{e^{\alpha_i}=1\}$ by $d\in \{0,2,4,\dots\}$. But now consider 
the operator $L(k)$ of example \ref{ex:L}, and take $\xi_1=
\frac{1}{2}\alpha_i^\vee|\alpha_i|$. It follows directly from the 
explicit definition \ref{dfn:duop} for $T_{\xi_1}$ that such a 
meromorphic function can be an eigenfunction of $L(k)$ only 
if 
\begin{equation}\label{eq:indwall}
d(d+1-2k_i)=0
\end{equation}
(In other words, the operator $L(k)$ has exponents $0$ and 
$(1-2k_i)/2$ (in the sense of Oshima \cite{Osh})) along the wall  
$\{e^{\alpha_i}=1\}$ (considered in the orbit space 
$W\backslash H^{reg}$). Hence for generic $k$ it is 
clear that we must have $d=0$. But an irreducible componentof the set of
singularities  of a meromorphic function cannot have codimension $> 1$, 
hence the result is true for arbitrary $k$.
\qed

\begin{remark}
The first part of Corollary \ref{cor:cont} is remarkable,
and it is not so easy to prove directly
for Harish-Chandra series without the deformation theory in $k$. 
The reason is that in the situation of a symmetric space, the 
two exponents of $L(k)$ along a wall are $0$ and a nonpositive integer 
(by \ref{eq:indwall}). In this case there possibly exist true meromorphic 
solutions, but by the deformation in $k$ it is clear that this 
possibility does not occur for the linear combination of Harish-Chandra 
series considered in the Corollary. 
\end{remark}


\subsection{The hypergeometric function} 
The function $\tilde{F}$ is more beautiful and well behaved than $\Phi$. When 
normalized at $e\in H$ this function will be denoted $F(\lambda,k;h)$, and 
this function will be called the hypergeometric function for the root 
system $R$. It is the natural generalization of the elementary spherical
function on a symmetric space with restricted root system $R$ (compare with 
Example~\ref{ex:L}).

\begin{theorem}{\rm (\cite{Op1}, Theorem 2.8)}
$\tilde F$ extends to an entire function
of $\lambda,\, k$ and $h$  (in a tubular neighbourhood of $A$).
\end{theorem}
\pf From Lemma \ref{lem:pole} and the explicit formula
for the $c$-function it is clear that $\tilde F$ may have 
first order poles along hyperplanes of the form 
$(\lambda,\alpha^\vee)=n$. First consider the case $n=0$. In this case the
first order pole has to be removable since $\tilde{F}$ is $W$ invariant in
$\lambda$. Next if $n\not= 0$ we may assume that $\alpha=\alpha_i$ is simple and 
$n>0$
by $W$ invariance. Take the residue $\mathrm{Res}_{n,i}$ of $\tilde F$  
at the hyperplane $H_{n,\alpha_i}$.
Clearly $\mathrm{Res}_{n,i}$ is also a solution of \ref{eqn:hgs}, defined on a 
tubular neigbourhood of $A$ in $H$ and $W$ invariant there. Let $W_i$ be the rank one parabolic
subgroup $W_i=\{1,r_i\}$ and let $W^i$ denote 
the set of elements $w$ in $W$ such that $l(wr_i)>l(w)$.
By 
\ref{lem:pole}, there exists an asymptotic expansion on $A_+$ of the form ($\lambda \in
H_{n,\alpha_i}$):
\[
\mathrm{Res}_{n,i}(a)=\sum_{w\in W}d_{w}(\lambda,k)\Phi(w\lambda,k,a)
\]
with $d_w=0$ if $w\in W^i$ (and in particular, $d_e=0$).
The remaining leading exponents have, for generic $\lambda\in H_{n,\alpha_i}$,
no mutual differences in $P$. Hence we may, for any $j\in\{1,\dots,n\}$, 
separate $\mathrm{Res}_{n,i}$ into subsums  \[
\mathrm{\Sigma}_{w,j}(a)=\sum_{x\in W_j}d_{xw}(\lambda,k)\Phi(xw\lambda,k,a)
\]
using the monodromy action of the $\theta_v$ (see text preceding 
Corollary \ref{cor:monoind}) such that $r_j v= v$.
By Lusztig's formula \ref{eqn:lus}
we have $[T_j,\theta_v]=1$ 
for such $v$. Hence these subsums $\mathrm{\Sigma}_{w,j}$ 
are still $T_j$ invariant. Therefore, the boundary value of $\mathrm{\Sigma}_{w,j}$
along the wall $e^{\alpha_j}=1$ is 
a multiple of an ordinary hypergeometric function. From the theory of
asymptotic expansion of the ordinary hypergeometric function we obtain that 
$d_w=d_{r_jw}=0$ if either $d_w=0$ or 
$d_{r_jw}=0$. This, combined with the prior remark that $d_e=0$, 
implies that $d_w=0\ \forall w\in W$, by a simple inductive argument on the 
length of $w$. Hence the
pole at $H_{n,\alpha_i}$ was removable. \qed 
\index{Gauss summation formula}
\begin{theorem}{\rm(Gauss summation formula \cite{Op3})}
\label{thm:gauss}
The function $\tilde F$ can be evaluated explicitly at the 
unit element of $H$:
 $\tilde{F}(\lambda,k;e) = \tilde{c}(\rho(k),k)$.
This evaluation is equivalent to the following limit 
formulae:
When $k_\alpha \le 0$ for all $\alpha$,
then
\[
\lim_{a \in A_+, a \downarrow e} \Phi(\lambda,k;a)
= \tilde{c}(-\lambda, 1-k).
\]
\end{theorem}
\pf 
We normalize 
\mbox{\index{f@$F(\lambda,k;a)$}}
\[
F(\lambda,k;a)
 = \frac{1}{\tilde{c}(\rho(k),k)} \tilde{F}(\lambda,k;a)
\]
and consider the value at the identity $f(\lambda,k):= F(\lambda,k;e)$.
It follows from Theorem \ref{thm:shiftprop}(d) that,
since
\[
G_-(k+1) \tilde{F}(\lambda,k+1) = \tilde{F}(\lambda,k),
\]
one has in any case the property that
$f(\lambda,k)$ is entire and periodic in $k$.
One can show $f(\lambda,k)$ is nonvanishing. 
We also see that $f(\lambda,k)\in {\bf R}$ if $\lambda,k$ are real.
Finally one can show that $k \mapsto f(\lambda,k)$ is entire with 
growth order $\le 1$. (This is technical, but essentially based on 
the recurrence relations (\ref{(**)}) for $\Gamma_\kappa(\lambda,k)$.)
By Hadamard's factorization theorem for entire functions 
one concludes that a function with these properties must be 
constant in k, and therefore $f(\lambda,k)=f(\lambda,0)=1$ for all 
$\lambda$ and $k$. For the formulation in terms of the limits of 
Harish-Chandra series: consult \cite{Op3}.
\qed 
\begin{dfn}
\index{hypergeometric function}
$F(\lambda,k;a)$ is called the {\em hypergeometric function}
for the root system $R$.
\end{dfn}
%
\newpage
\section{The KZ connection}

The goal of this section is to understand properly the analogue of the 
polynomials $E(\lambda,k)$ for arbitrary $\lambda\in\mathfrak h^*$. 
We call this analogue {\it nonsymmetric hypergeometric
functions}\index{nonsymmetric hypergeometric
functions}. The construction of nonsymmetric local solutions 
of the $T_\xi$ on a $W$-orbit leads naturally to the study of the 
so called Knizhnik-Zamolodchikov connection. We will gain a lot 
of insight in the equations \ref{eqn:hgs} by doing this exercise. Most 
importantly perhaps, it will become plain that the system has regular 
singularities. Also, it will naturally bring into play the action of the
affine  Weyl group by virtue of the affine intertwiners of Cherednik as
discussed  in Section 4.

Basic references for this section are \cite{Op} and \cite{Hec}. 

\subsection{Nonsymmetric hypergeometric functions}

For each element $h\in H^{\text{reg}}$, we define 
\index{sl@$S(\lambda,k)$}
\[
S_{Wh}(\lambda,k)
=\left\{\varphi\in{\cal{O}}_{Wh}\,;\, 
p(T_\xi(k))\varphi=p(\lambda)\varphi\enskip \text{\rm for any} \enskip p\in
S(\mathfrak h)^W\right\}. \]

\begin{prop} 
The space $S_{Wh}(\lambda,k)$ is an $\mathbf{H}(R_+,k)$-module and the
dimension of the subspace $S_{Wh}(\lambda,k)^W$ of $W$-invariant elements
is $|W|$. 
\end{prop}
\begin{pf}
Recall that $\mathbf{H}(R_+,k)$ is realized as the algebra generated
by $W$ and $\{T_\xi(k)\,;\, \xi\in\mathfrak h\}$ and also that the center of
$\mathbf{H}(R_+,k)$ is $\{p(T_\xi(k))\,;\, p\in S(\mathfrak h)^W\}$ (Lemma
4.1). 
Hence, $S_{Wh}(\lambda,k)$ is a module for $\mathbf{H}(R_+,k)$. 
By definition of $D_p$ (Definition \ref{def:de}), 
$S_{Wh}(\lambda,k)^W$ is the space 
of solutions of the hypergeometric system (\ref{eqn:hgs}). 
Then, by Theorem \ref{thm:rank}, $\dim S_{Wh}(\lambda,k)^W=|W|$. \qed
\end{pf} 

We now want to understand the weight subspace 
\index{sll@$S_{Wh}(\lambda,k)^\lambda$}
\[ 
S_{Wh}(\lambda,k)^\lambda
=\left\{\varphi\in{\cal{O}}_{W h}\,;\, 
T_\xi(k)\varphi=\lambda(\xi)\varphi\enskip\text{\rm for any}\enskip
\xi\in\mathfrak h\right\}.\] 
We have a map from $S_{Wh}(\lambda,k)^\lambda$ to $S_{Wh}(\lambda,k)^W$
given by $\varphi\mapsto\sum_{w\in W}\varphi^w$. 
(As in Section~3, we use the notation $\varphi^w=\varphi(w^{-1}\cdot)$ 
for a function $\varphi$).   
The following simple algebraic lemmata serve to prove that this is an
isomorphism if
$\lambda$ satisfies some conditions.

\begin{lemma} The $\mathbf{H}$-module 
$I_\lambda=\text{\rm Ind}_{S(\mathfrak h)}^{\mathbf{H}}(\C_\lambda)$ 
is called the minimal principal series module induced from the character
$\lambda$.  It is isomorphic to the regular representation as 
$\C[W]$-module.
Suppose that $\lambda$ satisfies $\lambda(\alpha^\vee)\not=0, \pm k_\alpha$ for
all $\alpha\in R_+$. Then $I_\lambda$ is the direct sum of its 
one dimensional weight spaces $I_\lambda^\mu$ with $\mu\in W\lambda$.   
Moreover, $I_\lambda$ is irreducible and the map
\[p:I_\lambda^\mu\ni v\mapsto \sum_{w\in W}wv\in I_\lambda^W\] 
is an isomorphism for any $\mu\in W\lambda$. Finally, every module 
over $\mathbf{H}$ with central character $\lambda$ and dimension $\leq |W|$ 
is ismorphic to $I_\lambda$. 
\end{lemma}
\begin{pf}
Under the assumption on $\lambda$ we see that the kernel of the 
intertwiners $I_w$ cannot have a nontrivial intersection with the 
weight space $I_\lambda^\lambda$. Hence all weight spaces of the 
form $I_\lambda^\mu$ with $\mu\in W\lambda$ are at least one dimensional.
Thus by a dimension count every weight space $I_\lambda^\mu$ is one
dimensional,  and the intertwiners $I_w$ act as isomorphisms. The
irreducibility  of $I_{\lambda}$ follows from the remark that any nonzero
submodule has to contain  at least one weight vector, but we have seen that
all weight vectors  are cyclic. Suppose that $0\not=v\in I_\lambda^\mu$ and
that $p(v)=0$. Then  $\mathbf{H}v=\C[W]v$ has dimension less than $|W|$,
contradicting the irreducibility. If $M$ is a module with central character 
$\lambda$ and dimension $\leq |W|$, then we argue as before that all its 
weight spaces with weight $\mu\in W\lambda$ have dimension 1. In particular, 
there is a nonzero weight vector of weight $\lambda$, which gives rise to 
an isomorphism with $I_\lambda$.
\qed 
\end{pf}

\begin{lem}\label{lem:Wlam}
Let $M$ be any $\mathbf{H}(R_+,k)$-module with central character
$\lambda$. Denote by $M^\lambda$ the weight space with weight
$\lambda$ and by $M^W$ the subspace of $W$-invariant elements. 
If $\lambda(\alpha^\vee)
\not=0,\pm k_\alpha$ for all $\alpha\in R_+$, then 
$M$ is semisimple and isotypic of type $I_\lambda$.
The map 
\[p:M^\lambda\ni v\mapsto \sum_{w\in W}wv\in M^W\] 
is an isomorphism. If $M^W$ is finite dimensional then 
$M$ itself is finite dimensional with 
$\mathrm{dim}(M)=|W|\mathrm{dim}(M^W)$.
\end{lem} 
\begin{pf} 
For a given $v\in M$ let us consider the submodule  
$\mathbf{H}v$. This is a quotient of the module
$Q_\lambda=\mathbf{H}/J_\lambda$ with $J_\lambda$ the ideal generated by the
central elements  $p-p(\lambda)$ with $p\in S(\mathfrak(h))^W$. 
It is clear that $Q_\lambda$ can be represented by $\mathfrak{H}\otimes
\C[W]$ with $\mathfrak{H}$ the harmonic elements in $S(\mathfrak(h))$.
Hence $Q_\lambda^W$ has dimension $|W|$, and for every $q\in Q_\lambda^W$,  
$Hq$ is isomorphic to $I_\lambda$ by the previous lemma. Thus $Q_\lambda$ is 
a direct sum of $|W|$ copies of $I_\lambda$. Now everything claimed follows
from the previous lemma.
\qed
\end{pf}

\begin{rem}\label{rem:inv}
 The inverse of 
\[p:M^\lambda\ni v\mapsto \sum_{w\in W}wv\in M^W\]
is given by the application of the element 
$q\in S(\mathfrak{h})$ given by \[q=\prod_{\alpha\in R_+}
\left(1-\frac{k_\alpha}{\lambda(\alpha^\vee)}\right)^{-1}
\prod_{w\in W, w\not=e}
\frac{\xi-w\lambda(\xi)}{\lambda(\xi)-w\lambda(\xi)},\]
where $\xi$ is any element in $\mathfrak h$ satisfying
$\lambda(\xi)\not=w\lambda(\xi)$ for all $w\not=e$
\end{rem}
\begin{pf}
It is sufficient to prove this for $M=I_\lambda$.
Consider the following identity in $\C[W]$:
\[
|W|\epsilon^+=\sum_wc_w(\lambda){\tilde I}_w(\lambda)
\]
(notations as in Remark \ref{rem:coc} and Lemma \ref{lm:equiv}).
We compute the coefficients $c_w$ easily by the following 
remarks. First of all, one verifies directly that 
\[
c_{w_0}(\lambda)=\prod_{\alpha\in
R_+}\frac{\lambda(\alpha^\vee)+k_\alpha}{\lambda(\alpha^\vee)}.
\]
Using the cocycle relation of Remark \ref{rem:coc} and the observation 
$\epsilon^+\cdot {\tilde I}_w(\lambda)=\epsilon^+$ it
follows  that $c_w(\lambda)=c_{w_0}(w_{0}w\lambda)$, hence 
 \[
c_w(\lambda)=\prod_{\alpha\in
R_+}\frac{w\lambda(\alpha^\vee)-k_\alpha}{w\lambda(\alpha^\vee)}
\] 
Apply this decomposition of $p=|W|\epsilon^+$ to $v=1\in I_\lambda^\lambda$
and we see that $q\circ p(1)=1$, as desired.
\qed \end{pf}

\begin{cor} Retain the assumptions of Lemma \ref{lem:Wlam}.
The dimension of $S_{Wh}(\lambda,k)$ is $|W|^2$, and this defines
a local system $S(\lambda,k)$ of $\mathbf{H}=\mathbf{H}(R_+,k)$-modules with
central character  $\lambda$ on the regular orbit space.   
The monodromy of this local system centralizes the $\mathbf{H}$-module
structure, and gives $S_{Wh}(\lambda,k)$ the structure of a
$\mathbf{H}^\mathrm{aff}(R_+,q)$-module. More precisely, $S_{Wh}(\lambda,k)$ is the
direct sum of $|W|$ copies  of the monodromy of the equations \ref{eqn:hgs}.
\end{cor} 
\begin{pf}
We leave to the reader the easy verification that monodromy of $S(\lambda,k)$ 
commutes with the actions on $S(\lambda,k)$ by $W$ and by Dunkl operators.
By the previous lemmata, $S(\lambda,k)$ is the direct sum of weight
spaces $S(\lambda,k)^\mu$ all of which are isomorphic to  $S(\lambda,k)^W$ 
via the intertwiner $p$ for the monodromy. (And of course, 
$S(\lambda,k)^W$ is nothing but the local system of solutions of
\ref{eqn:hgs}).
 \qed
\end{pf}
\begin{cor}\label{cor:nonsymhgf} 
If $\Re k_\alpha\geq0$ for any $\alpha\in R_+$, then there exists a
unique holomorphic function 
$G(\lambda,k;\cdot)$ \index{g@$G(\lambda,k;\cdot)$} in a tubular
neighbourhood of $A$ such that 
\begin{align*}
&T_\xi(k)G(\lambda,k;\cdot)=\lambda(\xi)G(\lambda,k;\cdot),\tag1\\
&G(\lambda,k;e)=1.\tag2
\end{align*}
\end{cor}

\begin{pf} 
For $\lambda$ satisfying $\lambda(\alpha^\vee)\not=0, \pm k_\alpha$ 
for any $\alpha\in R_+$, we define 
\[G(\lambda,k;\cdot)=|W|D_qF(\lambda,k;\cdot).\]
By Remark \ref{rem:inv}, (1) is clear. 

Since this function satisfies (again by Remark \ref{rem:inv}):
\[F(\lambda,k;\cdot)=
\frac1{|W|}\sum_{w\in W}G^w(\lambda,k;\cdot),
\] 
(2) follows from Theorem \ref{thm:gauss}. 
The apparent poles in $\lambda$ are
removable  because of the next lemma, from which the uniqueness also follows. 
\end{pf}

\begin{lem}\label{lem:nul}
Let $\varphi\in S(\lambda,k)^\lambda$ be a holomorphic function 
in a neighbourhood 
of $e\in A$. If $\Re k_\alpha\geq0$ for any $\alpha\in R_+$, 
then $\varphi(e)=0$ implies $\varphi=0$. 
\end{lem}
\begin{pf} 
Let $\{\xi_i\}$ be an orthonormal basis of $\mathfrak a$ and let 
$\{\xi_i^*\}$ be the dual basis. The lowest homogeneous part of the 
operator 
\[\sum_{i=1}^n\xi_i^*T_{\xi_i}(k)=\sum_{i=1}^n\xi_i^*\partial_{\xi_i}
+\sum_{\alpha\in R_+}\frac{k_\alpha\alpha}{1-e^{-\alpha}}
(1-r_\alpha)\]
at the origin is equal to 
\[E(k)=\sum_{i=1}^n\xi_i^*\partial_{\xi_i}
+\sum_{\alpha\in R_+}k_\alpha(1-r_\alpha).\]
Assume that $\varphi\not=0$ and let $f$ be the lowest homogeneous part 
of $\varphi$ with degree $m\geq0$. 
By the equation $\sum_{i=1}^n\xi_i^*T_{\xi_i}(k)\varphi=\lambda\varphi$, 
we have $E(k)f=\left(m+\sum_{\alpha\in R_+}k_\alpha(1-r_\alpha)\right)f
=0$. Since $\C[W]f$ is a $\C[W]$-module, 
we can express $f$ as a sum $\sum_{\delta\in\hat{W}}f_\delta$ 
of $\delta$-equivariant parts $f_\delta$ for each $\delta\in\hat{W}$. 
The element $\sum_{\alpha\in R_+}k_\alpha(1-r_\alpha)$ is central in 
$\C[W]$, hence acts on an irreducible $\C[W]$-module $\delta$ 
by a scalar. 
It is easy to see that this scalar is equal to 
\[\epsilon_\delta(k)=\sum_{\alpha\in R_+}
k_\alpha(1-\chi_\delta(r_\alpha)/\chi_\delta(e)), \]
where $\chi_\delta$ is the character of $\delta$, and we have the 
following equation: 
\[(m+\epsilon_\delta(k))f_\delta=0\quad\text{\rm for each}
\enskip \delta\in\hat{W}.\]
On the other hand, since $\Re\epsilon_\delta(k)$ is not less than zero 
for each $\delta\in\hat{W}$ by assumption, 
we have $f_\delta=0$ unless $m=0$. Contradiction. \qed
\end{pf}

We shall prove the removability of poles of $G(\lambda,k)$. 
Assume that $G(\lambda,k)$ has a singularity. 
Since $F(\lambda,k)$ is an entire function of $(\lambda,k)$ and 
by the expression $G(\lambda,k)=|W|D_qF(\lambda,k)$, 
$G(\lambda,k)$ is meromorphic in $(\lambda,k)$ and its singular set 
is the zero set of a function that depends only on $(\lambda,k)$. 
Let $(\lambda_0,k_0)$ be a regular point and let $\varphi$ be 
an irreducible holomorphic function in a neighbourhood $V$ of 
$(\lambda_0,k_0)$ such that the zero set of $\varphi$ is equal to 
the singular set in $V$. 
Let $l\in\N$ be the smallest integer such that 
$\tilde{G}=\varphi^lG$ extends holomorphically to $V$. 
By continuity and the property (2), $\tilde{G}(\lambda,k,e)=0$ 
for any singular point $(\lambda,k)$ in $V$ and, by Lemma \ref{lem:nul}, 
$\tilde{G}(\lambda,k)\equiv0$ for these points. This is a contradiction. 
\qed

\begin{ex}
Let us consider the $\mbox{BC}_1$ case, i.e. $R=\{\pm\alpha, \pm2\alpha\}.$ 
We use the notation in Example~5.3. 
The functions $F$ and $G$ are expressed as follows: 
\[\begin{cases}
F(\lambda,k;x)={}_2F_1(a,b,c;z),\\
G(\lambda,k;z)={}_2F_1(a,b,c;z)+\frac1{4b}(y-y^{-1})
{}_2F_1'(a,b,c;z),\\
\end{cases}\]
where, ${}_2F_1(a,b,c;z)$ is  Gauss' hypergeometric function. 
\index{Gauss' hypergeometric function}  
\end{ex}

\begin{rem} 
We have seen that $p: S^\lambda\rightarrow S^W$ is an 
isomorphism if $\lambda(\alpha^\vee)\not=0,\pm k_\alpha$ for all 
$\alpha\in R$,  
and that this map is an intertwiner for the monodromy 
representation of $\mathbf{H}^\mathrm{aff}(R_+,q_i)$. 
In fact, for sufficiently generic parameters, 
we have two isomorphisms: 
\begin{align*}
S(\lambda,k)&\simeq I_\lambda^{|W|}
\quad(\text{\rm as $\mathbf{H}$-module}),\\
&\simeq\left(\text{\rm Ind}_{\C[Q^\vee]}^{\mathbf{H}^\mathrm{aff}}
e^{2\pi\sqrt{-1}(\lambda)}\right)^{|W|}
\quad(\text{\rm as $\mathbf{H}^\mathrm{aff}$-module}).
\end{align*}
These two actions commute with each other. 
Notice that also the shift operators 
$G_\pm(k): S(\lambda,k)^W\rightarrow S(\lambda,k\pm1)^W$ 
and the intertwiners 
$I_w: S(\lambda,k)^\lambda\rightarrow S(w\lambda,k)^{w\lambda}$ 
($w\in W^e$) commute with the 
$\mathbf{H}^\mathrm{aff}$-action. 
\end{rem}

\begin{rem}\label{rem:act}
Since $T_\xi(k)$ is {\it not} $W$-equivariant, 
$G(w\lambda,k;a)$ and $G^w(\lambda,k,a)$ do not coincide. 
The correct relationship between them is given by affine 
intertwiners: 
\[I_wG(\lambda,k)
=\left(\prod_{a\in R_+^a\cap w^{-1}R_+^a}
(\lambda(a)+k_a)\right)G(w\lambda,k)\quad
\text{\rm for }\enskip w\in W^e.\]
\end{rem}

\subsection{The role of the Knizhnik-Zamolodchikov connection} 

Let \index{omegal@$\Omega^l$} $\Omega^l$ be the sheaf 
of holomorphic $l$-forms on 
$\mathfrak h^{\text{\rm reg}}$. 
We use the notation $\Omega_h^l$ and $\Omega_{Wh}^l$ 
analogously to ${\cal{O}}_h$ and ${\cal{O}}_{Wh}$. 

Define an operator 
\index{dl@$d(\lambda,k)$} 
$d(\lambda,k): \Omega_{Wh}^l\rightarrow\Omega_{Wh}^{l+1}$ 
by 
\[d(\lambda,k)=d-d(\lambda+\rho(k))
+\sum_{\alpha\in R_+}k_\alpha(1-e^{-\alpha})^{-1}d\alpha\otimes(1-r_\alpha).\]

As in the proof of Lemma \ref{lem:nul}, let $\{\xi_i\}$ be an orthonormal
basis  of $\mathfrak a$ and let $\{\xi_i^*\}$ be its dual basis of $\mathfrak
a^*$.  Since the action of $d(\lambda,k)$ is expressed as 
\begin{align*}
d(\lambda,k)&(\varphi\otimes dx_1\wedge\dots\wedge dx_l)\\
=&\sum_{i=1}^n\left(\partial_{\xi_i}-(\lambda+\rho(k))(\xi_i)
+\sum_{\alpha\in R_+}\frac{k_\alpha\alpha(\xi_i)}{1-e^{-\alpha}}
(1-r_\alpha)\right)\varphi\\
&\qquad\quad\otimes d\xi_i^*\wedge dx_1\wedge\cdots
\wedge dx_l\\
=&\sum_{i=1}^n\left(T_{\xi_i}(k)-\lambda(\xi_i)\right)
\varphi\otimes d\xi_i^*\wedge dx_1\wedge\cdots\wedge dx_l,\\
\end{align*}
we have $d(\lambda,k)^2=0$, and 
\[0\longrightarrow S^\lambda\overset{\text{\rm inj.}}\longrightarrow
{\cal{O}}_{Wh}\overset{d(\lambda,k)}\longrightarrow
\Omega_{Wh}^1\overset{d(\lambda,k)}\longrightarrow
\Omega_{Wh}^2\longrightarrow\dots\]
is a cochain complex. 

Note that $\Omega_{Wh}^l$ is isomorphic to 
$(\Omega_{Wh}^l\otimes\C[W])^W$ by 
\[\Omega_{Wh}^l\ni\varphi \overset{\sim}\longmapsto
\sum_{w\in W}\varphi^w\otimes w\in(\Omega_{Wh}^l\otimes\C[W])^W.
\]
On the other hand, $(\Omega_{Wh}^l\otimes\C[W])^W$ 
is also isomorphic to $\Omega_h^l\otimes\C[W]$ by 
\[\Omega_h^l\otimes\C[W]\ni\varphi\otimes v
\overset{\sim}\longmapsto
\sum_{w\in W}\varphi^w\otimes wv\in
(\Omega_{Wh}^l\otimes\C[W])^W.\]
Via these isomorphisms, we have a new cochain complex: 
\[0\longrightarrow {\cal{L}}^\lambda\overset{\text{\rm inj.}}\longrightarrow
{\cal{O}}_{h}\otimes\C[W]\overset{\nabla(\lambda,k)}\longrightarrow
\Omega_{h}^1\otimes\C[W]\overset{\nabla(\lambda,k)}\longrightarrow
\Omega_{h}^2\otimes\C[W]\longrightarrow\dots\]
Since the isomorphism $\Omega_{Wh}^l\overset{\sim}\longrightarrow
\Omega_h^l\otimes\C[W]$ is given by 
\[(\varphi_w)_{w\in W}\mapsto\sum_{w\in W}\varphi^w_{w^{-1}}\otimes w\quad
(\varphi_w\in\Omega_{w\cdot h}^l)\] 
and the inverse is given by 
\[\sum_{w\in W}\psi_w\otimes w\mapsto(\psi_{w^{-1}}^w)_{w\in W},\]
the operator \index{nabla@$\nabla(\lambda,k)$} 
$\nabla(\lambda,k)$ is expressed as follows: 
\begin{align*}
\nabla(\lambda,k)(\psi\otimes w\otimes &  dx_1\wedge\dots\wedge dx_l)\\
=&\sum_{i=1}^n\nabla_{\xi_i}(\lambda,k)(\psi\otimes w)
\otimes d\xi_i^*\wedge dx_1\wedge\dots\wedge dx_l,
\end{align*}
with
\begin{align*}
\nabla_\xi(\lambda,k)=&w\left(T_{w^{-1}\xi}(k)-w\lambda(\xi)\right)w^{-1}
\quad\text{\rm (multiplication in $\mathbf{H}(R_+,k)$)}\\
=&\partial_\xi+\frac12\sum_{\alpha\in R_+}
k_\alpha\left(\alpha(\xi)\frac{1+e^{-\alpha}}{1-e^{-\alpha}}
\otimes(1-r_\alpha)+\alpha(\xi)\otimes r_\alpha\epsilon_\alpha\right)
-w\lambda(\xi),
\end{align*}
and
$\epsilon_\alpha(w)=-\text{\rm sgn}(w^{-1}\alpha)w$. 
The last expression is a consequence of (\ref{eqn:wxi}), 
and the reflections in $\nabla_\xi(\lambda,k)$ act on $\C[W]$ 
by left multiplication. 


\begin{defn} 
\index{Knizhnik-Zamolodchikov connection}
\index{KZ-connection}
We call the coinvariant derivative $\nabla(\lambda,k)$ 
the (trigonometric) {\em Knizhnik-Zamolodchikov connection} 
(KZ-connection in the sequel).
\end{defn}

\begin{cor}{\rm (Matsuo \cite{M})}\label{cor:matsuo} 
The KZ connection is integrable and the map 
$\sum_{w\in W}\psi_w\otimes w\mapsto\sum_{w\in W}\psi_w$ 
gives an isomorphism from ${\cal{L}}_\lambda$ to $S^W$ 
if $\lambda(\alpha^\vee)\not=0,\pm k_\alpha$ for any $\alpha\in R$. 
\end{cor}
The isomorphism in Corollary~\ref{cor:matsuo} is called 
the Matsuo isomorphism. \index{Matsuo isomorphism}

\begin{rem} 
We can easily extend this isomorphism to the weaker condition 
``$\lambda(\alpha^\vee)\not=k_\alpha$ for any $\alpha\in R_+$''. 
\end{rem}

\begin{rem} 
By Corollary~\ref{cor:nonsymhgf}, 
$G(\lambda,k)\in S(\lambda,k)^\lambda$. 
Then, by the above discussion, 
the vector $\sum_{w\in W}G^w(\lambda,k)\otimes w$ is an element 
of 
\index{llambda@${\cal{L}}_\lambda$}
\[
{\cal{L}}_\lambda
=\{\psi\in{\cal{O}}_h\otimes\C[W]\,;\, 
\nabla(\lambda,k)\psi=0\}.\]
\end{rem}
\newpage
%
\section{Harmonic Analysis on $A$}
In this section we study the eigenfunction transform 
$\F$  for the algebra of Dunkl operators acting on $C_c^\infty(A)$.  
We shall prove a Paley-Wiener theorem and an explicit inversion 
formula for $\F$, when $k_\alpha\geq 0$ for all $\alpha\in R$. 
The transform $\F$ was called the Cherednik transform 
in \cite{Op} and the Opdam transform in \cite{C1}. 
We will simply use the generic name ``Fourier transform'' here. 
\subsection{Paley-Wiener theorem}
For $f,\,g\in C_c^\infty(A)$, define
\index{(@$(f,g)_k$}
\[(f,g)_k=\int_A f(a)\overline{g(a)}\delta_k(a)da,\]
where
\index{deltak@$\delta_k(a)$}
\[\delta_k(a)=\prod_{\alpha\in R_+}
\left\vert a^{\alpha/2}-a^{-\alpha/2}\right\vert^{2k_\alpha}\]
and $da$ is the Lebesgue measure on $A$ normalized by 
$\mbox{vol}(A/\exp(Q^\vee))=1$. 
In this section we assume that $k_\alpha\geq 0$ for 
all $\alpha\in R$. In this and the next section we shall only give 
complete proofs when there is something new to add to the 
ideas in the literature. Otherwise we shall content ourselves 
with references. 

The following lemma is an easy computation.
\begin{lem}{\rm (\cite[Lemma 7.8]{Op})}
\[(T_\xi f,g)_k=(f,(-w_0 T_{w_0(\bar{\xi})}w_0)g)_k.\]
Here $w_0$ is the longest element in $W$. 
\end{lem}
\begin{definition}For $f\in C_c^\infty(A)$ and $\lambda\in\h^*$, define
\[\F(f)(\lambda)=\int_A f(a)G(-w_0\lambda,k;w_0 a)\delta_k(a)da.\]
And for $\varphi$ a ``nice function'' on $\h^*$, define
\index{j@$\J$}
\[\J(\varphi)(a)=\int_{\i\a^*}
\varphi(\lambda)G(\lambda,k;a)\sigma(\lambda)d\mu(\lambda),\]
where
\index{sigma@$\sigma(\lambda)$}
\[\sigma(\lambda)=\prod_{\alpha\in R_+}
\frac{\Gamma(\lambda({\alpha}^\vee)+k_\alpha)
\Gamma(-\lambda({\alpha}^\vee)+k_\alpha+1)}
{\Gamma(\lambda({\alpha}^\vee))\Gamma(-\lambda({\alpha}^\vee)+1)},\]
and $d\mu(\lambda)$ is the translation invariant holomorphic $n$-form 
such that the volume of $\i\a/2\pi \i P$ equals $1$. 
\end{definition}

First we need to show that $C_c^\infty(A)$ is mapped by 
\index{f@$\F$} $\F$ in 
a space of nice functions, so that the composition 
$\J\F(f)$ makes sense.
\index{c@$C_a$} 
Given $a\in A$, let $C_a$ denote the convex hull of $Wa$ and 
\index{ha@$H_a$} let $H_a$ denote the support function given by 
\[H_a(\lambda)=\mbox{sup}\{\lambda(\log b)\,;\,b\in C_a\}.\]
\index{Paley-Wiener type}
An entire function $\varphi$ on $\h^*$ is said to be of 
{\em Paley-Wiener type} $a$ if
\[{}^\forall N\in \N,\,{}^\exists C>0\,:\,
|\varphi(\lambda)|\leq C(1+|\lambda|)^{-N}\exp({H_a(-\mbox{Re}(\lambda))})
\quad (\lambda\in\h^*).\]
Let \index{pwa@$PW(a)$}$PW(a)$ be the space of entire functions of 
the Paley-Wiener type
$a$ and \index{pw@$PW$} $PW=\bigcup_{a\in A}PW(a)$. 

\begin{thm}{\rm (\cite[Proposition 6.1, Corollary 6.2]{Op})}
\label{thm:estimate}
For all $k\in K^{\text{reg}}$ (here regular means:
$\tilde{c}(\rho(k),k)\not=0$) and all compact subset $D$ of $A$,  and all
$p\in S(\h)$, there exists $C>0$ and $N\in \N$ such that  \[\sup_{a\in
D}|\partial(p)G(\lambda,k;a)| \leq C(1+|\lambda|^N)\exp({\max_w\{\mbox{\rm
Re}(w\lambda(\log a))\}}).\] \end{thm}
\pf (Sketch) 
If $a$ and $\xi$ are regular elements in the same Weyl chamber, 
we can see that
\[\partial_\xi(a^{-2\mu}\sum_w|G(\lambda,k,w^{-1}a)|^2)\leq 0\]
from KZ connection, where $\mu\in W\mbox{Re}\lambda$ such that 
$\mu(\xi)=\max_w\{\mbox{\rm Re}(w\lambda(\xi))\}$. 
This proves the theorem for $p\equiv 1$. 
The statement for  general $p\in S(\h)$ follows from Cauchy's formula.
\qed

\index{Paley-Wiener theorem}
\begin{thm} {\rm (Paley-Wiener theorem \cite[Theorem 8.6]{Op})}
\par\noindent
{\rm (a)} $\F\,:\,C_c^\infty(C_a)\rightarrow PW(a)$
\par\noindent
{\rm (b)} $\J\,:\,PW(a)\rightarrow C_c^\infty(C_a)$
\end{thm}
\pf 
(a) follows directly from Theorem \ref{thm:estimate}. 
\noindent
Using asymptotic expansion (b) can be proved in the same way as 
Helgason's proof of the 
Paley-Wiener theorem for Riemannian symmetric spaces \cite{Hel}.
\qed
\subsection{Inversion and Plancherel formula}
\begin{thm}{\rm (see \cite{Op})}\label{thm:inv}
$\F\J$ and $\J\F$ are identical on $PW$ and $C_c^\infty(A)$ respectively.
\end{thm}
\begin{pf}
The theorem was first proved by Opdam\cite{Op}. 
Here we will give an outline of Cherednik's 
proof of Theorem \ref{thm:inv} (\cite{C1}). It is a very nice proof, 
based on the action of the affine intertwiners. The nonsymmetric theory 
is essential now.

One checks by direct computation that 
\begin{eqnarray}
\F(I_i f)(\lambda) & = & -(\lambda(a_i)+k_i)\F(f)(r_i\lambda), 
\label{eqn:71} \\
\F(T_\xi f)(\lambda) & = & \lambda(\xi)\F(f)(\lambda),
\label{eqn:72}
\end{eqnarray}
Combined these formulae show that 
\begin{equation}\label{eqn:73}
\F(f^{r_i})=\F(f)^{r_i}-k_i\frac{\F(f)-\F(f)^{r_i}}{a_i}=Q_i(\F(f)).
\end{equation}
Here \index{qi@$Q_i$} $Q_i$  
is the Lusztig operator, which is the action of $r_i$ in 
the module
\[\mbox{Ind}_{\C [W]}
^{PW\otimes_{S(\h)}\H(R_+,k)}(\mbox{triv}).\]

Next one checks that 
\begin{equation}\label{eqn:74}
\J(Q_i(\varphi))=\J(\varphi)^{r_i}\quad i=0,1,\dots,n.
\end{equation}
This is delicate if $i=0$, since we need a contour shift here (the
proof for $i\not=0$ is the same, but without the shift).  
If  $i=0$ it is only
true for $k_\alpha\geq  0$ ($\alpha\in R$). For the proof we need
\begin{equation}\label{eqn:75}
\left(1+\frac{k_i}{\lambda(a_i)}\right)\sigma(\lambda)=
\left(1-\frac{k_i}{\lambda(a_i)}\right)\sigma(r_i\lambda),
\end{equation}
which follows easily from the definition of $\sigma$. 

We have
\begin{gather*}
\begin{split}
\J&(Q_i(\varphi))
 = \int_{\i\a^*} Q_i(\varphi)(\lambda)
G(\lambda,k;a)\sigma(\lambda)d\mu(\lambda) \\
  &=  \int_{\i\a^*} \left(\varphi^{r_i}-
k_i\frac{\varphi-\varphi^{r_i}}{\lambda(a_i)}\right)
G(\lambda,k;a)\sigma(\lambda)d\mu(\lambda) \quad 
(\mbox{by }(\ref{eqn:73}))  \\
  &= \int_{\i\a^*} \varphi({r_i}\lambda)
\left(1+\frac{k_i}{\lambda(a_i)}\right)
G(\lambda,k;a)\sigma(\lambda)d\mu(\lambda) \\
& -k_i\int_{\i\a^*} \varphi(\lambda)\frac{1}{\lambda(a_i)}
G(\lambda,k;a)\sigma(\lambda)d\mu(\lambda) \\
 &= \int_{r_i(\i\a^*)} \varphi(\lambda)
\left(1+\frac{k_i}{\lambda(a_i)}\right)
G(r_i\lambda,k;a)\sigma(\lambda)d\mu(\lambda) \\
& -k_i\int_{\i\a^*} \varphi(\lambda)\frac{1}{\lambda(a_i)}
G(\lambda,k;a)\sigma(\lambda)d\mu(\lambda) 
\quad(\mbox{by }(\ref{eqn:75})) \\
 &= \int_{\i\a^*} \varphi(\lambda)
\left(\left(
1+\frac{k_i}{\lambda(a_i)}\right)G(r_i\lambda,k;a)
-\frac{k_i}{\lambda(a_i)}G(\lambda,k;a)\right)
\sigma(\lambda)d\mu(\lambda)  \\
 &= \J(\varphi)^{r_i}.\\
\end{split}
\end{gather*}

In last steps we use shift of contour for $i=0$ and 
a formula for $G^{r_i}$ based on the formula for $I_iG$ (cf. Remark
\ref{rem:act}):  \[G^{r_i}(\lambda,k;a)=\left(
1+\frac{k_i}{\lambda(a_i)}\right)G(r_i\lambda,k;a)
-\frac{k_i}{\lambda(a_i)}G(\lambda,k;a).\]
Observe that the necessary shift of contour when $i=0$ is 
allowed when $k_\alpha>0$, since the only pole of $\sigma$ 
that possibly needs to be reckoned with is cancelled by the 
factor
\[1+\frac{k_0}{\lambda(a_0)}=
\frac{1-\lambda(\theta^\vee)+k_\theta}{1-\lambda(\theta^\vee)}. 
\]
However, when $k_\alpha<0$ the poles at $\lambda(a_i)+k_i$ enter 
into the positive chamber, and these destroy the symmetry for $i=0$.

By (\ref{eqn:73}) and (\ref{eqn:74}), 
$\J\circ\F$ commutes with action of $W^e$ on $C_c^\infty(A)$. 
In particular, $\J\circ\F$ commutes with multiplications 
by $e^\lambda\,(\lambda\in P)$. It is easy to see that the 
ideal $i_{x_0}$ of functions in $C_c^\infty(A)$ that vanish at some
point $x_0\in A$ can be written as $j_{x_0}C_c^\infty(A)$, where
$j_{x_0}$ denotes the maximal ideal at $x_0$ in $\C[P]$. Hence $\J\circ\F$
maps $i_{x_0}$ into itself, for all $x_0$. Therefore it has to be
multiplication by a $f\in C^\infty(A)$. Since $\J\circ\F$ is also $W$
equivariant, $f$ must be $W$ invariant. Finally, by (\ref{eqn:72}), it has to
also commute with $T_\xi$-action on $C_c^\infty(A)$. Thus we have \[T_\xi
f=\partial_\xi f=0\quad \mbox{for all } \xi,\] and $f$ must be a constant. 
One can prove that the constant is $1$ by considering the asymptotics. 

Conversely $\F\circ \J$ commutes with multiplications by 
polynomials $p\in S(\h)$. As before, $\F\circ\J$ has to be 
multiplication by some function $g$. Computing $\J\F\J(\varphi)$ in 
two ways, we have 
\[\J(\varphi)=\J(g\varphi).\]
At $e\in A$ we have
\[\int_{\i\a^*}\varphi(\lambda)\sigma(\lambda)d\mu(\lambda)
=\int_{\i\a^*}g(\lambda)\varphi(\lambda)\sigma(\lambda)d\mu(\lambda),\]
hence $g\equiv 1$. 
\qed
\end{pf}
The inversion formula we have derived now is NOT the inversion formula of 
the spectral decomposition of $\C_c^\infty(A)$ for the action of the
commutative algebra of Dunkl-Cherednik operators (this algebra of operators 
is not even closed with respect to the $*$ operator!). Accordingly, the 
function $\sigma$ is not positive (not even real), we have no 
Plancherel formula and no extension of $\F$ to an $L_2$ space.  
One can fix this by considering the decomposition of $\C_c^\infty(A)$ 
with respect to its structure as a pre-unitary module of the 
action of the noncommutative $*$ algebra $\mathbf{H}$, and this 
point of view was used in \cite{Op}. A simpler way out of this is the 
reduction of the
transform to the $|W|$-symmetric situation.  If $f\in C_c^\infty(A)$ is
$W$-invariant, then \begin{equation}\label{eqn:symf}
\F(f)(\lambda)=\int_A f(a)F(-\lambda,k;a)da,
\end{equation}
which coincides with the \index{Harish-Chandra transform} 
Harish-Chandra transform 
for spherical functions if the parameter $k$ corresponds 
to the root multiplicities of a Riemannian symmetric space. 

The $W$-invariance of $f$ results in the $W$-invariance 
of $\F(f)$. Replacing $G$ by $F$ in the transform $\J$, we have
\begin{equation}\label{eqn:syminv}
f(a)=\int_{\i\a^*}\F(f)(\lambda)F(\lambda,k;a)\sigma'(\lambda)
d\mu(\lambda),
\end{equation}
where
\index{sigmap@$\sigma'(\lambda)$}
\[
\sigma'(\lambda)
=\prod_{\alpha\in R_+}
\frac{\Gamma(\lambda({\alpha}^\vee)+k_\alpha)
\Gamma(-\lambda({\alpha}^\vee)+k_\alpha)}
{\Gamma(\lambda({\alpha}^\vee))\Gamma(-\lambda({\alpha}^\vee)).}\]
Notice that 
\[\sigma'(\lambda)=\frac{1}{c(\lambda,k)c(-\lambda,k)}=\frac{1}{|c(\lambda,k)|^2}>0,\]
where
\index{c@$c(\lambda,k)$}
\[
c(\lambda,k)
=\frac{\tilde{c}(\lambda,k)}{\tilde{c}(\rho(k),k)}.\]
Formula (\ref{eqn:syminv}) is a $k$-deformation of Harish-Chandra's  
inversion formula for spherical transform. 
For arbitrary $k$ 
($k_\alpha\geq 0,\,\alpha\in R$) it had been conjectured by Heckman
and Opdam and was proved by Opdam\cite{Op}. For group case, see 
\cite[Ch IV]{Hel}. 
%
\newpage
%
\section{The attractive case (Residue Calculus)}
In the previous section we gave the inversion formula for $\F$ 
for the repulsive case, $k_\alpha\geq 0$ for all $\alpha\in R$. 
In this section we consider the attractive case, $k_\alpha<0$ for 
all $\alpha\in R$ (cf. \cite{Op2}). The spectral decomposition 
involves lower dimensional spectra.

\subsection{Paley-Wiener theorem and Plancherel theorem}
The formula
\[(f,g)_k=\int_A f(a)\overline{g(a)}\delta_k(a)da,\]
gives an inner product only as long as $\delta_k(a)$ is 
locally integrable.

\begin{thm}{\rm ( \cite[Proposition 5.1]{Hec}, 
\cite[Proposition 1.1]{Op2})} \label{thm:k}
$\delta_k(a)$ is locally integrable if and only if $k$ is in the 
connected component of $\{k\,;\,\tilde{c}(\rho(k),k)>0\}$ containing
$k_\alpha\geq 0$ for all $\alpha\in R$. 
In particular this is satisfied in the following two situations:
\begin{itemize}
\item[{(a)}] $k_\alpha\geq 0$ for all $\alpha\in R$.
\item[{(b)}] $k_\alpha<0$ for all $\alpha\in R$, and 
$\rho(k)({\theta}^\vee)+k_\theta+1>0$.  
\end{itemize}

Here, as always, $\theta$ is the highest short root.
In case (a), $\delta_k(a)$ is locally integrable and in case (b),
$\delta_k(a)$ is even integrable.
\end{thm}

\begin{rem}
If $R$ is  simply laced, the condition for $k$ in the theorem 
means that $k>-1/d_n$, where $d_n$ is the Coxeter number. 
\end{rem} 
\begin{rem}
If $\delta_k(a)$ is integrable, then $G(-\rho(k),k,\cdot)=1$ is square
integrable  with respect to $\delta_k(a)da$. On the other hand, in the sense 
of the previous section its Fourier transform is zero. Clearly 
the inversion formula with purely continuous spectrum as in the 
previous section now fails! 
\end{rem}

From now on we assume that we are in the situation of Theorem 
\ref{thm:k}(b) (the so-called attractive case). And we will restrict 
ourselves to the $W$-symmetric case, in view of the remarks made in the 
last part of the previous section.

We define 
$\F$ as before, but we define $\J$ by 
\begin{equation}\label{eqn:syminv2}
(\J\varphi)(a)=\int_{\gamma+\i\a^*}\varphi(\lambda)
\Phi(\lambda,k;a)\frac{d\mu}{c(-\lambda,k)},\quad \varphi\in PW,
\end{equation}
where 
\index{asm@$\a_-^*$}
$\gamma\in\a_-^*=\{\lambda\in
\a^*\,;\,\lambda(\alpha^\vee)<0\,{}^\forall
\alpha>0\}$
 such that $\gamma({\alpha}^\vee)<k_\alpha$  
and $a\in A_+$. 
By Lemma \ref{lem:pole}, $\Phi(\lambda,k;a)$ is holomorphic in $\lambda$
if $\mbox{Re}(\lambda({\alpha}^\vee))<1-\varepsilon$ for all
$\alpha\in R_+$ and $\varepsilon>0$. 
If $k_\alpha\geq 0$ for all $\alpha\in R$, then  
(\ref{eqn:syminv2}) coincides with 
the right hand side of 
(\ref{eqn:syminv}) for $\F f=\varphi$ by analytic continuation 
and symmetrization.  

As we have seen, the proof of Theorem \ref{thm:inv} by Cherednik fails. 
However, the original proof of the inversion formula survives:
\begin{thm} {\rm (see \cite[Theorem 5.4]{Op2})} Still $\J\F$ and $\F\J$ are 
identical.
\end{thm}

We will now engage a process to refine the defining formula for $\J$ in
such a way that $\J$ becomes integration of $\lambda$ over some subset of
$\mathfrak{h}^*$, against the  kernel  $F(\lambda,k;a)$ multiplied by a 
positive measure, the Plancherel measure. This will give rise to the extension
of $\F$ to  $L_2(A,\delta_k)^W$, and eventually to an isometric
isomorphism of $L_2(A,\delta_k)^W$ with the $L_2$ space on
$\mathfrak{h}^*$ defined by the Plancherel measure. In other words, this leads
to the spectral resolution of the commutative algebra of differential
operators $D_p$, $p\in S(\mathfrak{h})^W$. 

This process consists of a shift of the
contour of (\ref{eqn:syminv2}) from $\gamma+\i\a^*$ to $\i\a^*$.  
The residual contours one encounters along the way also move as though they
are attracted by the origin (and these again pick up residues along the way,
and so on). When everybody comes to a standstill, we have contours of
integration in  every possible dimension. Next we have to symmetrize, and then
finally we will have the integral defining $\J$ satisfying the properties
described mentioned above.


Let us first formulate the results of all this precisely.
We need some
terminology:  
\begin{definition}\label{def:res}
An affine subspace \index{l@$L$} $L\subset\a^*$ 
is called {\em residual} if
\index{residual}
\begin{equation}
\#\{\alpha\in R\,;\,{\alpha}^\vee(L)=k_\alpha\}
=\#\{\alpha\in R\,;\,{\alpha}^\vee(L)=0\}+\mbox{codim}(L).
\end{equation}
Notice that $\a^*$ itself is residual. 
If a residual subspace $L$ is a point, we call it a 
{\em distinguished point}. \index{distinguished point}
 Given $L$ residual, let \index{cl@$c_L$} $c_L$ denote 
the orthogonal projection 
of $0\in\a^*$ on $L$, and put
\[L=c_L+V^L,\]
\index{ltemp@$L^{\text{temp}}$}
\[L^{\text{temp}}=c_L+\i V^L\subset\h^*.\]
\end{definition}
\begin{rem}
The classification of residual subspaces reduces to the classification of 
distinguished points by ``parabolic
induction''. If $k_\alpha=k$ for all
$\alpha\in R$, the distinguished points correspond to the distinguished
nilpotent orbits in the semisimple Lie algebra $\g_\C({R}^\vee)$. 
Such orbits were classified by Carter and Bala.
\end{rem}

The desired formula for $\J$ is given in the next theorem:
\begin{thm}{\rm (\cite[Theorem 3.4]{Op2})}\label{thm:res}
\[
\J\varphi(a)=\sum_{L}
\int_{L^{\text{\rm temp}}}\varphi(\lambda)F(\lambda,k;a)
d\nu_L(\lambda,k).
\]
Here 
\index{dnu@$d\nu_L(\lambda,k)$}
\begin{equation}\label{eqn:meas}
d\nu_L(\lambda,k)
=\gamma_L(k)f_L(\lambda,k)d\mu_L(\lambda),
\end{equation}
\[f_L(\lambda,k)=\tilde{c}(\rho(k),k)^2 
\frac{\prod'_L \Gamma(\lambda({\alpha}^\vee)+k_\alpha)}
{\prod'_L \Gamma(\lambda({\alpha}^\vee))},
\]
$\mu_L$ is Lebesgue measure on $L^{\text{\rm temp}}$ such 
that $\text{\rm vol}(\i V^L/2\pi\i(P\cap V^L))=1$,  
$\prod_L'$ is the product of the $\Gamma$-factors of the 
roots which do not vanish identically on $L$,  
$0\leq \gamma_L(k)\in{\bf Q}$, and the sum is taken over 
all the residual subspaces $L$ such that $c_L\in \a_-^*$. 
\end{thm}
\begin{cor}{\rm (\cite[Theorem 5.7, Corollary 5.8]{Op2})}
\label{cor:sqint}
$\nu_L(\lambda,k)$ is a 
positive measure (if nonzero).
The $W$-invariant square integrable eigenfunctions of $L(k)$ 
are $F(\lambda(k),k;\cdot)$ with $\lambda(k)$ distinguished in 
$\a_-^*$ and 
$\gamma_{\lambda(k)}(k)>0$. For these we have
\begin{gather*}
\begin{split}
\int_A F(\lambda(k),k;a)^2\delta(k,a)da & \\
=\pm\gamma_L^{-1}|W\lambda(k)|^{-1} & 
\frac{\prod_{\alpha\in R_+}\Gamma(\rho(k)({\alpha}^\vee)+k_\alpha)^2
\prod_{\alpha\in R\setminus R_z}\Gamma(\lambda(k)({\alpha}^\vee))}
{\prod_{\alpha\in R_+}\Gamma(\rho(k)({\alpha}^\vee))^2
\prod_{\alpha\in R\setminus
  R_p}\Gamma(\lambda(k)({\alpha}^\vee)+k_\alpha)}, 
\end{split}
\end{gather*}
where 
\begin{eqnarray*}
R_z & = &\{\alpha\in R\,;\,\lambda(k)({\alpha}^\vee)=0\text{ for all
  }k\}, \\
R_p & = &\{\alpha\in R\,;\,\lambda(k)({\alpha}^\vee)+k_\alpha=0 \}. \\
\end{eqnarray*}
\end{cor}
The parameters $\lambda(k)$ in Corollary \ref{cor:sqint} 
are classified in \cite[Section 4]{HO}. 

\begin{ex}\rm (see \cite{BHO})
If $k_\alpha=k$ for all $\alpha\in R$ then  
$\lambda(k)=\rho(k)$ is distinguished and for  $F(\rho(k),k;\cdot)=1$, 
we have
\[
\int_A \delta(k,a)da=
\prod_{i=1}^n\left(\begin{array}{c}d_ik \\ k\end{array}\right)
\frac{\pi}{\sin(-m_i\pi k)},\]
where $m_i$ are the exponents and $d_i=m_i+1$ are the degrees.
\end{ex}

In the rest of the section, we 
will give an outline of the proof of Theorem \ref{thm:res}. 

\subsection{Residues}
Given a finite arrangement of affine hyperplanes $\mathcal{H}$ 
in a Euclidean space $V$, we choose for each $H\in\mathcal{H}$ a vector
$\alpha_H\in V$, and a number $k_H\in\R$ such that 
\[H=\{\lambda\in V\,;\,(\alpha_H,\lambda)=k_H\}.\]
Let $\L$ be the lattice of intersections of elements of $\mathcal{H}$, 
ordered by inclusion (and $V\in L$ by definition). Let $\omega$ be a
rational $n$-form on $V_\C$, with poles possibly at the hyperplanes of
$\mathcal{H}$, but nowhere else. Let $PW$ denote the space of Paley-Wiener
functions, with rapid decay in the imaginary direction.

\bigskip\noindent
\underline{GOAL} Study the functional
\index{xv@$X_{V,\gamma}$}
\[X_{V,\gamma}\,:\,PW\rightarrow\C, \,
\varphi\mapsto\int_{\gamma+\i V}\varphi\omega,\]
in particular what happens when $\gamma$ moves from chamber 
to chamber.

\bigskip
We may rewrite $X_{V,\gamma}$ in many different ways as a sum of 
$X_{V,\gamma^\prime}$'s and residual integrations over lower dimensional
contours. In fact, we will describe a systematic way of pointing out a special
chamber in each $L\in\mathcal{L}$, 
to which we want to move $\gamma$. The point
is that this defines a unique way of rewriting $X_{V,\gamma}$. 

Given $L\in \L$, let $c_L$ be the orthogonal projection of $O\in V$
onto $L$. Write $L=c_L+V^L$, where $V^L\subset V$ is a linear
subspace and $\mathcal{C}=\{c_L\,;\,L\in \L\}$, the set of 
centers. The next lemma is elementary, but very effective.
\index{center}

\begin{lem}{\rm (\cite[Lemma 3.1]{Op2})}
There exists a unique collection of tempered distributions on $X_c$,
$c\in \mathcal{C}$ 
such that \\
{\rm
 (a)} $\mbox{\rm supp}(X_c)\subset\cup_{L;\,c_L=c}\i V^L$, \\
{\rm (b)} $X_c$ has finite order, \\
{\rm (c)} $X_{V,\gamma}(\varphi)=\sum_{c\in \mathcal{C}}
X_c(\varphi(c+\cdot))$ for all $\varphi\in PW$.
\end{lem}
The distributions $X_c$ play a crucial role. We refer to $X_c$ as 
``the local contribution of $X=X_{V,\gamma}$ at the center $c$''.
\begin{rem}
The value of $X_c$ does not change when either $0$ or $\gamma$ passes 
a hyperplane that does not contain $c$. Hence, 
when computing $X_c$, we may always assume that both $O$ and $\gamma$ 
are in chambers which contain $c$ in their closure. In other words, 
we reduce in this way to consider the central arrangement of hyperplanes
that contain the center $c$.
\end{rem}

\begin{lem}{\rm (\cite[Lemma 3.3]{Op2})}\label{lem:sup}
\index{h@$\mathcal{H}$}
Let $\mathcal{H}$ be a central arrangement with center $c$. If $X_c\not=0$, 
then $O$ must be in the closure of the antidual chamber of the
chamber in which $\gamma$ lies. Explicitly,
\[O\in\overline{\sum_{H\in\mathcal{H}'}\R_+c_H+\sum_{H\in\mathcal{H}''}\R_-c_H}+c,\]
where $\mathcal{H}'$ 
is the set of non-separating hyperplanes for $\gamma$ and
$O$, and $\mathcal{H}''=\mathcal{H}\setminus \mathcal{H}'$.
\end{lem}
The above result follows from the next example, the special case of 
normal crossings, since every arrangement of hyperplanes can be approximated 
by arrangements with normal crossings only. In this normal crossing case it 
is a simple exercise using the geometry of simplicial cones.
\begin{ex} {\rm (normal crossing case)} Suppose $(\gamma,\alpha_H)<k_H$ for all
$H\in \mathcal{H}$, and $\mathcal{H}$ is divisor with normal crossings at
$c=\cap_{H\in\mathcal{H}}H$. Assume 
\[\omega=\prod_{H}((\lambda,\alpha_H)-k_H)^{-1}d\lambda\]
and assume that $O$ is in the antidual of $\gamma$. Then 
\begin{gather*}
\begin{split}
X_c(\varphi(c+\cdot))
&=(-2\pi\i)^n\det(\alpha_H,\alpha_{H'})^{-1/2}\varphi(c) \\
&=(-2\pi\i)^n\frac{1}{\text{vol}(V/\sum_H\Z\alpha_H)}\varphi(c).
\end{split}
\end{gather*}
\end{ex}

\subsection{The arrangement of shifted root hyperplanes}
Assume that we have a root system $R$, irreducible, reduced, in
$V=\a^*$, and root multiplicities $k_\alpha\in\R_-$. Let
$R^\vee\subset\a$ be the set of coroots, and normalize the Lebesgue
measure $dx$ (resp. $d\lambda$) on $\a$ (resp. $\i\a^*$) 
such that $\mbox{covol}(Q^\vee)=1$ (resp.  $\mbox{covol}(2\pi\i
P)=1$).  Denote by \index{cp@$c'(\lambda,k)$} 
$c'(\lambda,k)$ the rational function
\[c'(\lambda,k)=\prod_{\alpha\in R_+}
\frac{\lambda({\alpha}^\vee)+k_\alpha}{\lambda({\alpha}^\vee)}.\]

Consider 
\index{yv@$Y_{{\mathfrak{a}^*},\gamma}(\varphi)$}
\begin{gather*}
\begin{split}
X_{{\mathfrak{a}^*},\gamma}(\varphi)& =\int_{\gamma+\i\a^*}\varphi(\lambda)
\frac{d\lambda}{c'(-\lambda,k)}, \\
Y_{{\mathfrak{a}^*},\gamma}(\varphi)
& =\int_{\gamma+\i\a^*}\varphi(\lambda)
\frac{d\lambda}{c'(-\lambda,k)c'(\lambda,k)}, 
\end{split}
\end{gather*}
where $\gamma\in\a^*_-$ such that $\gamma({\alpha}^\vee)-k_\alpha<0$
for all $\alpha\in R_+$. Let 
\[H_\alpha=\{\lambda\in\a^*\,;\,\lambda({\alpha}^\vee)=k_\alpha\,
{}^\forall\alpha\in R\}\]
and let $\mathcal{C}$ be the set of centers of the corresponding 
intersection lattice $\L$. For $c\in\mathcal{C}$, denote by $X_c$ and 
$Y_c$ the local contribution of $X_{{\mathfrak{a}^*},\gamma}$ and $Y_{{\mathfrak{a}^*},\gamma}$. 
Given $c\in\mathcal{C}$, denote by $W_c$ the stabilizer in $W$ of $c$,
and let $A_c$ denote the symmetrization operation
\index{ac@$A_c\varphi(\lambda)$}
\[
A_c\varphi(\lambda)
=|W_c|^{-1}\sum_{w\in W_c}c'(w\lambda,k)
\varphi(w\lambda).\]
Notice that this is holomorphic in a neighbourhood of $c+\i {\mathfrak{a}^*}$ if
$\varphi$ is so. 

\begin{lem}{\rm (\cite[Proposition 3.6]{Op2})}
For $c\in\mathcal{C}\cap\overline{\a^*_-}$ and $w\in W$, we have 
\[X_{wc}=Y_c\circ w^{-1}\circ A_{wc}.\]
\end{lem}
This has the following application, which is of substance when $c$ is
singular.  Suppose that $\lambda$ is in the support  of some $Y_c$ with $c\in
\mathfrak{a}^*_-$. If $w\lambda$ is not in the support of $X_{wc}$ then
$A_{wc}\varphi(w\lambda)$ must be zero.  By  Lemma \ref{lem:sup} this is
always the case when $wc=Re(w\lambda)\not\in\overline{{}_-\mathfrak{a^*}}$.
This argument will show that the hypergeometric function $F(\lambda,k,\cdot)$
has all  its leading exponents in $\overline{{}_-\mathfrak{a^*}}$ for such
$\lambda$, hence is tempered by a well known criterion of Casselman and
Milici\'c. This is the content of Corollary \ref{cor:temp}. Let us now
formulate this argument on a technical level.
The next result is a direct application of Lemma \ref{lem:sup}.
\begin{cor}
{\rm(\cite[Corollary 3.7]{Op2})}\label{cor:av} For $c\in\mathcal{C}$, write 
\[{}_-{\mathfrak{a}^*}^c= \sum_{\alpha\in
R_+,c({\alpha}^\vee)=k_\alpha}\R_-\alpha\subset 
\overline{{}_-{\mathfrak{a}^*}},\] where $\overline{{}_-{\mathfrak{a}^*}}$ is 
the closure of antidual of ${\mathfrak{a}^*}_+$.  Let $c\in
\mathcal{C}\cap\overline{{\mathfrak{a}^*}_-}$ and $w\in W$ with
$wc\notin{}_-{\mathfrak{a}^*}^{wc}$. If $\lambda\in c+\mbox{\rm supp}(Y_c)$
then $A_{wc}\varphi(w\lambda)=0$ for all $\varphi\in PW({\mathfrak{a}^*}_c)$. 
\end{cor}

First of all, recall that in this attractive case $k_\alpha<0$, we are
interested only in the situation where $\delta_k(a)$ is integrable on
$A$, and we have seen that this means that condition (2) in
Theorem~\ref{thm:k} holds. It means  geometrically that 
\[C_{\rho(k)}\subset\{\lambda\in\a^*\,;\,
|\lambda({\alpha}^\vee)|<1+k_\alpha\,{}^\forall\alpha\in R\}.\]
Choose an open convex $W$-invariant set $U$ between these sets.

\begin{lem}{\rm (\cite[Proposition 2.2]{Op2})}
Let $a\in A_+$. Then $\lambda\mapsto\Phi(\lambda,k;a)$ is holomorphic
on $\a_-+U+\i \mathfrak{a}^*$, and uniformly bounded there.
\end{lem}

\begin{lem}{\rm (\cite[Lemma 3.3]{Op2})}
Write $c(\lambda,k)=c'(\lambda,k)c''(\lambda,k)$. Then
$c''(\lambda,k)^{\pm}$ are holomorphic on $U+\i a^*$, and 
$c''(\lambda,k)^{-1}$ bounded, $c''(\lambda,k)$ of moderate
growth. Also $c''(-\lambda,k)^{-1}$ is holomorphic in $\a_-+U+\i\a^*$
and $c''(\lambda,k)c''(-\lambda,k)$ and $c'(\lambda,k)c'(-\lambda,k)$ 
are $W$-invariant. \end{lem}

\begin{lem}{\rm (\cite[Lemma 3.2]{Op2})}
All centers $c\in\mathcal{C}$ lie in $U$. 
\end{lem}

Corollary \ref{cor:av} contains important information about the
hypergeometric function, because the operator $A_{wc}$ plays a role in
its definition. If $c=Re(\lambda)$ then 
\begin{gather*}
\begin{split}
{F}(\lambda,k;a)
&=\sum_{w \in W} {c}(w\lambda,k) \Phi(w\lambda,k;a)\\
&=\sum_{w\in W/W_c}|W_c|A_{wc}(c''(\cdot,k)\Phi(\cdot,k;a))(w\lambda).
\end{split}
\end{gather*}
Together with the above results concerning the good behaviour 
of $\Phi$ and $c''$ on $U+i\mathfrak{a}^*$ this finally leads to 
the desired result:

\begin{cor}{\rm (\cite[Corollary 3.7]{Op2})}\label{cor:temp}
If $\lambda\in c+\mbox{\rm supp}(Y_c)$, $c\in \mathcal{C}\cap\a^*_-$,
and $w\in W$ such that $wc\notin{}_-{\mathfrak{a}^*}^{wc}$, then $a\mapsto
F(\lambda,k;a)$ is tempered on $A$. If $L=c$, and $Y_c\not=0$, then
$F(c,k;a)$ has exponential decay; such $F$ are called cuspidal.
\end{cor}

Now we need to say more about the shifted root hyperplane
arrangement. There are two very special geometric 
peculiarities of this arrangement that make everything work 
properly. It is obvious that the local contributions of
$Y_{\mathfrak{a}^*,\gamma}$ have support at subspaces that are 
residual in the following sense.

\begin{definition}\label{def:res>}
$L$ is called residual if \index{residual}
\[
\#\{\alpha\in R\,;\,{\alpha}^\vee(L)=k_\alpha\}
\geq\#\{\alpha\in R\,;\,{\alpha}^\vee(L)=0\}+\mbox{codim}(L).
\]
\end{definition}
However, as we have seen in Definition \ref{def:res}, whenever the 
above inequality holds it has to be an equality! This is of crucial 
importance because this shows that the local contributions of 
$Y_{\mathfrak{a}^*,\gamma}$ are in fact densities (distributions 
of order 0). Another important point is that a residual subspace $L$ 
of dimension $k$ is determined by a distinguished point of 
a parabolic subsystems of rank $n-k$. In fact $L^{\mathrm{temp}}$ is the 
space of the corresponding unitary parabolic induction parameters, as embedded
in the parameter space of the minimal principal series. This stucture makes it
possible to work with ``unitary parabolic inducion''. 
The second peculiarity has to do with the positivity of the relative
Plancherel measures on $L^{\mathrm{temp}}$ needed in this inductive process. Here one
needs the property that $-c_L$ and $c_L$ are in the same orbit of the fixator
group of ${\mathfrak{a}^*}^L$ in $W$.

The following theorem is proved by
the classification (!) of distinguished points.  
\begin{thm}\label{thm:clas}{\rm (\cite[Theorem 3.9,
Theorem 3.10, Remark 3.11]{HO})} If $L$ is residual in the sense of Definition
\ref{def:res>}, then  
\[
\#\{\alpha\in R\,;\,{\alpha}^\vee(L)=k_\alpha\}
=\#\{\alpha\in R\,;\,{\alpha}^\vee(L)=0\}+\mbox{\rm codim}(L).
\]
If $L$ is residual, then its center $c_L$ is a distinguished point for
$R_L=\{\alpha\in R\,;\,L({\alpha}^\vee)=constant\}$ and $-c_L\in
W(R_L)c_L$. 
\end{thm}
As indicated, this leads to:
\begin{cor}
If $L$ is residual, $c_L\in \mathcal{C}\cap \a^*_-$ and
$Y_{c_L}\not=0$, then it is in fact a measure, namely integration over
$c_L+\i {\mathfrak{a}^*}^L$ against the density 
\[d\nu_L'(\lambda,k)=\gamma_L(k)
\frac{\prod'|\lambda({\alpha}^\vee)|}
{\prod'|\lambda({\alpha}^\vee)+k_\alpha|}
d_L(\lambda),\]
where $\prod'$ denotes the product over all $\alpha\in R$, omitting
zero factors.
\end{cor}
The Corollary \ref{cor:temp} makes it possible to show (by induction, 
starting with the distinguished points) that all densities involved are 
in fact positive measures (and Theorem \ref{thm:clas} is crucially needed in
the inductive process):
\begin{cor}
The function $(c''(\lambda,k)c''(-\lambda,k))^{-1}$ is positive, 
bounded and real analytic on $c_L+\i {\mathfrak{a}^*}^L$, and 
$\nu_L(\lambda,k)=(c''(\lambda,k)c''(-\lambda,k))^{-1}\nu'_L(\lambda,k)$ 
is given by formula (\ref{eqn:meas}). It is a positive, real analytic 
measure when $\gamma_L(k)\not=0)$. 
\end{cor}

\begin{cor}
If $\varphi$ is a $W$-invariant, PW-function and
$\gamma({\alpha}^\vee)<k_\alpha$ for all $\alpha\in R_+$, then 
\begin{gather*}
\begin{split}
\J_\gamma(\varphi)&=X_{\a^*,\gamma}(\varphi\Phi(\cdot,k,a)
c''(-\lambda,k)^{-1}) \\
&=\int_{\gamma+\i\a^*}\varphi(\lambda)\Phi(\lambda,k,a)
\frac{d\lambda}{c(-\lambda,k)} \\
&=\sum_{L:\text{\rm residual},c_L\in\mathcal{C}\cap\a_-^*}
\int_{L^{\mathrm{temp}}}\varphi(\lambda)F(\lambda,k;a)d\nu_L(\lambda,k)\\
&:=\J(\varphi).
\end{split}
\end{gather*}
\end{cor}
Theorem \ref{thm:res} follows from this corollary. 
We finish with the main result, the Plancherel Theorem.  
\begin{theorem}{\rm (\cite[Theorem 5.5]{Op2})}
A residual subspace $L$ is called spherically tempered when 
$\nu_L\not=0$.
The map $\F$ extends 
naturally to an isometric isomorphism
\[
\F:L_2(A,\delta_k da)^W\to\left\{\bigoplus_{L\text{sph.temp.}}
L_2(L^{\mathrm{temp}},\nu_L(k))\right\}^W, 
\]
with inverse $\J$ as in Theorem 
\ref{thm:res}.
\end{theorem}
%

\input opdamam.ind
\end{document}